\documentclass{article}
\usepackage[%
journal=XXX,    
lang=british,   
]{ems-journal}

\usepackage{tikz}
\usetikzlibrary{matrix,positioning,decorations.pathreplacing}
\usepackage{pdflscape}
\usepackage{booktabs}
\usepackage{listings}
\usepackage{xcolor}
\usepackage{optidef}
\usepackage{algorithm}
\usepackage{algorithmic}
\usepackage{tcolorbox}
\tcbuselibrary{skins,breakable}

\newcommand{\rankk}{\ensuremath{\mathrm{rank}}}

\newcommand{\bM}{\mathbf{M}}

\newcommand{\cC}{\mathcal{C}}

\newcommand{\nrank}{\ensuremath{\mathrm{nrank}}}
\newcommand{\dist}{\ensuremath{\mathrm{dist}}}
\newcommand{\grad}{\ensuremath{\mathrm{grad}}}

\newcommand{\scalemath}[2]{\scalebox{#1}{\mbox{\ensuremath{\displaystyle #2}}}}

\DeclareMathOperator*{\argmin}{arg\,min}

\numberwithin{equation}{section}

\newtheorem{lemma}{Lemma}
\newtheorem{example}{Example}
\newtheorem{theorem}{Theorem}
\newtheorem{definition}{Definition}
\newtheorem{remark}{Remark}

\begin{document}

\title{Generic factorizations
and alternating direction method for fixed-rank approximation of matrix polynomials}
\titlemark{Short title for the running head - only if full title really does not fit (note the huge space)}


%


\emsauthor{1}{
    \givenname{Sweta}
    \surname{Das}
    \mrid{}
    \zblid{}
    \orcid{}}{S.~Das}

\emsauthor{2}{
    \givenname{Andrii}
    \surname{Dmytryshyn}
    \mrid{}
    \zblid{}
    \orcid{0000-0001-9110-6182}}{A.~Dmytryshyn}

\emsauthor{3}{
    \givenname{Froil\'an M.}
    \surname{Dopico}
    \mrid{}
    \zblid{}
    \orcid{0000-0001-7779-8960}}{F.~M.~Dopico}

\emsauthor{4}{
    \givenname{M{\aa}rten}
    \surname{Gulliksson}
    \mrid{}
    \zblid{}
    \orcid{}}{M.~Gulliksson}


\Emsaffil{1}{
    \department{School of Science and Technology}
    \organisation{\"{O}rebro University}
    \rorid{}
    \address{}
    \zip{701 82}
    \city{\"{O}rebro}
    \country{Sweden}
    \affemail{sweta.das@oru.se}}

\Emsaffil{2}{
    \department{1}{Department of Mathematical Sciences}
    \organisation{1}{Chalmers University of Technology
        and University of Gothenburg}
    \rorid{1}{}
    \address{1}{Chalmers Tv\"argata 3}
    \zip{1}{412 96}
    \city{1}{Gothenburg}
    \country{1}{Sweden}
    \department{2}{Corresponding author}
    \affemail{2}{andrii@chalmers.se}}

\Emsaffil{3}{
    \department{Departamento de Matem\'aticas}
    \organisation{Universidad Carlos III de Madrid}
    \rorid{03ths8210}
    \address{Avenida de la Universidad, 30
        (edificio Sabatini)}
    \zip{28911}
    \city{Legan\'es (Madrid)}
    \country{Spain}
    \affemail{dopico@math.uc3m.es}}

\Emsaffil{4}{
    \department{School of Science and Technology}
    \organisation{\"{O}rebro University}
    \rorid{}
    \address{}
    \zip{701 82}
    \city{\"{O}rebro}
    \country{Sweden}
    \affemail{marten.gulliksson@oru.se}}


%
%

\classification[65K10, 15A22, 15A54]{65F55}


\begin{abstract}
We present an algorithm for approximating a given matrix polynomial of degree $d$ by a matrix polynomial of prescribed rank and degree at most $d$. The method builds on recent advances in the theory of generic eigenstructures and factorizations of matrix polynomials with bounded rank and degree. We also give a specialized treatment of the case of matrix pencils (degree one matrix polynomials), motivated by its particular theoretical and practical importance. Lastly, numerical experiments are reported to assess the performance of the proposed approach and to compare it with existing methods.
\end{abstract}

\maketitle


\section{Introduction}

Low-rank approximation is one of the fundamental mechanisms by which large-scale linear-algebraic models are made tractable.  For constant matrices, the theory is classical: the singular value decomposition gives both a structural description
of low-rank matrices and, in unitarily invariant norms, optimal approximants.
This theory underlies a wide range of methods in numerical linear algebra,
systems theory, model reduction, system identification, and signal processing.

In many applications, however, the objects of interest are not constant but polynomial matrices, also known as matrix polynomials,  \cite{BHMS13, MaMT15, TiMe01}   $$P(\lambda)=\sum_{i=0}^{d} P_i\lambda^i,
        \qquad P_i\in \mathbb{C}^{m\times n}, \qquad i=0,\dots , d,$$
where ${\mathbb C}^{m\times n}$ is the set of $m\times n$ matrices with entries from the field of complex numbers $\mathbb{C}$. 

Matrix polynomials arise naturally in polynomial eigenvalue problems, linear systems and control, signal processing, and various forms of structured modeling \cite{BHMS13, TiMe01}. Their rank is the \emph{normal rank}, that is the rank of $P(\lambda)$ over the field of rational functions $\mathbb{C}(\lambda)$, and is denoted by $\nrank(P(\lambda))$. Equivalently,
$\nrank (P(\lambda)) = \max_{\mu \in \mathbb{C}} \rankk (P(\mu))$. The normal rank is a global property of a matrix polynomial and depends on all its coefficients. Consequently, the geometry of bounded-rank matrix polynomials is far more intricate than that of bounded-rank constant matrices, and this fundamentally changes the nature of the low-rank approximation problem. Despite the importance of matrix polynomials in applications, a reliable framework for their low-rank approximation has therefore remained elusive.

The purpose of this paper is to develop such a framework.  Our approach is based on recent advances in the theory of generic eigenstructures and factorizations of matrix polynomials of bounded normal rank \cite{DmDo17,DmDV23}.  These results identify generic families of matrix polynomials with prescribed rank and provide explicit parametrizations of such families.  We use these parametrizations to reformulate the constrained low-rank approximation problem
\begin{equation}
    \label{pro:minpol} 
            \min_{S(\lambda)}\ \|P(\lambda)-S(\lambda)\|_F, \quad \mbox{such that}
        \quad
        \deg (S(\lambda)) \le d,\quad \nrank (S(\lambda)) \le r,
\end{equation}
as an unconstrained problem (presented in \eqref{eq.unsconpoly1}), that allows efficient and reliable solution. 

Our framework applies to matrix polynomials of arbitrary degree and to any prescribed target rank $r<\min(m,n)$.  Its main feature is that the rank and degree constraints are built into the parametrization. Thus every polynomial produced by the algorithm has normal rank at most $r$, avoiding a common difficulty of losing the rank condition in methods that impose low rank only indirectly or through penalty terms. The resulting optimization problem is then approached using the alternating least-squares method. The parametrized formulation also reflects the intrinsic dimension of the rank-constrained family: as the target rank decreases, so does the number of parameters.  This makes the approach particularly attractive for computing very-low-rank approximations of matrix polynomials.

We deliberately refer to the construction as a framework rather than as a single algorithm.  It consists of several largely independent components, e.g., the choice of generic rank-$r$ families, the choice of a parametrization of these families, and the numerical method used to solve the resulting optimization problem.  Each of them can be modified, providing us with the flexibility to create algorithms tailored to particular applications, including structure-preserving approximations and more scalable optimization schemes.

The paper is organized as follows.  Section~\ref{sec:sota} reviews the current state of the art on low-rank approximation of matrix polynomials and related special cases.  In Sections~\ref{sec:pencase} and~\ref{sec:pol}, we develop the
proposed framework first for matrix pencils and then for matrix polynomials of arbitrary degree, including the necessary theoretical background.  Section \ref{sec:conv} discusses convergence properties of the resulting iterative schemes.  Numerical experiments illustrating the behavior and computational performance of the proposed approach are presented in Section~\ref{sec:numex}. We end with conclusions and some future research directions in Section~\ref{sec:conclusions}. 

\section{State-of-the-art for approximations with matrix polynomials}
\label{sec:sota}

We naturally divide this section by the degree of the considered polynomials. Note also that the majority of the results mentioned here are dealing with the partial case of the problem \eqref{pro:minpol}, defined by $r=\min(m,n)-1$.
For $m=n$, this case is known as the ``distance-to-singularity problem'' and is motivated by the need for robust computations; see, for example, \cite{GnGu25}, where the proximity to singularity causes instability in a delay differential equation.

{\bf The case of matrices, $d=0$.} 
The problem of approximating a given matrix $A$ by a matrix of prescribed rank $r<\min(m,n)$ is among the best understood problems in numerical linear algebra. The Eckart--Young--Mirsky theorem characterizes the optimal rank-$r$ approximation via the truncated singular value decomposition (SVD) \cite{GoVa13},  providing both the existence and the conditions of the uniqueness of the best approximant as well as a constructive algorithm.

\textbf{The case of matrix pencils, $d=1$.} 
The problem already becomes complicated in this setting. It was first introduced by Byers, He, and Mehrmann \cite{ByHM98} in the special case of $r=\min(m,n)-1$ and $m=n$. After the characterizations and upper and lower bounds for matrix pencils developed in \cite{ByHM98}, essentially no further progress was made for almost 20 years due to the intrinsic complexity of the problem \cite{KrVo15}. 
More recently, interest in this type of problem has been growing. In \cite{GuLM17}, Guglielmi, Lubich, and Mehrmann propose an approach for computing nearest singular pencils based on solving ordinary differential equations. Dopico, Noferini, and Nyman in \cite{DoNN23} suggest a method based on Riemannian optimization. It is important to emphasize that both \cite{DoNN23,GuLM17} address only the case $r=\min(m,n)-1$ and $m=n$ too. Several other works, which are also able to treat matrix polynomials, are discussed below.

\textbf{The case of matrix polynomials, $d\geq1$.} 
Naturally, low-rank approximation for matrix polynomials is significantly more challenging. In particular, extending the methods developed for the pencil case may be extremely difficult or even impossible. 
The work of Giesbrecht, Haraldson, and Labahn \cite{GiHL17} was the first to address the problem for matrix polynomials. More recently, several new approaches have been proposed in \cite{DaBo23, GnGu23, GnGu25, GNNP24}, leading to substantial progress in this area. The first named, introduced by Das and Bora \cite{DaBo23}, suggests methods based on perturbation theory. In \cite{GnGu25}, Gnazzo and Guglielmi formulate the distance to singularity problem ($r=\min(m,n)-1$) for general nonlinear matrix-valued functions as a nonconvex optimization problem with singularity imposed as a constraint; the polynomial case is described in detail in \cite{GnGu23}. The third approach, known as the Riemann-Oracle \cite{GNNP24}, developed by Gnazzo, Noferini, Nyman, and Poloni, provides a general-purpose Riemannian optimization framework for a broad class of matrix nearness problems, including the distance to singularity problem and, in Section 9, a strategy as to how to approach the case of any rank $r$, though an implementation is missing. 

%

\section{Approximation with low-rank matrix pencils}
\label{sec:pencase}

Matrix pencils are matrix polynomials of degree at most one. So, the approach presented in Section~\ref{sec:pol} for matrix polynomials of arbitrary degree covers the pencil case. Nevertheless, we provide a complete and self-contained treatment of this case in this section. The main motivation is the relative simplicity of matrix pencils, as well as their intrinsic importance. In fact, in the literature on distance problems and low-rank approximations, many contributions focus exclusively on the pencil case, see, e.g., \cite{ByHM98,DoNN23,GuLM17,KrVo15}.  We believe that presenting this case separately significantly simplifies the understanding of the general polynomial case, since the main steps of the method remain the same, but the technicalities are substantially reduced.

It is also worth noting that the set of bounded-rank matrix pencils admits a description in terms of factorizations (or in terms of a sum of pencils of rank at most one) that does not require taking closures of the involved subsets, in contrast to the general polynomial case (compare Theorems \ref{lem:P_as_union_of_C's} and \ref{thm.relationgeneigen}). Although this distinction is conceptually important, it does not play an essential role in the developments of our method.

In Sections \ref{sec:form_min_pen} and \ref{sec:sol_min_pen}, we formulate the minimization problems for matrix pencils using Theorem \ref{lem:P_as_union_of_C's} and solve them using the alternating least squares method, see, e.g., \cite{LiLR16, SrJa03}. The convergence properties of this approach are analyzed in Section \ref{sec:conv} in the general polynomial case. Before we proceed, some preliminary information about matrix pencils is summarized in the next subsection.

\subsection{Preliminary information about matrix pencils} \label{sec.prelimpenc}

All the matrices in this paper have complex entries, unless otherwise specified. To avoid repetitions, we define some basic concepts immediately for matrix polynomials. $\mathbb{C}[\lambda]$ denotes the ring of polynomials in one variable $\lambda$ with coefficients in $\mathbb{C}$, $\mathbb{C}[\lambda]^{m \times n}$ denotes the set of $m\times n$ matrix polynomials, i.e., the set of $m\times n$ matrices with entries from $\mathbb{C}[\lambda]$, and $\mathbb{C}[\lambda]^m$ stands for  $\mathbb{C}[\lambda]^{m \times 1}$, i.e., for the set of column vector polynomials with $m$ entries. The degree of a matrix polynomial $P(\lambda) \in \mathbb{C}[\lambda]^{m \times n}$ is the largest degree of its entries and is denoted by $\deg(P)$.
The Frobenius norm of $P(\lambda) = \sum_{i= 0}^d \lambda^i P_i$ is defined in terms of its matrix coefficients as follows $\|P(\lambda)\|_F = \left(\sum_{i= 0}^d \|P_i\|_F^2 \right)^{\frac{1}{2}}$. 


We denote by $\mathcal{P}^{m\times n}_r$ the set of $m\times n$ matrix pencils with normal rank at most $r$, i.e.,  
$\mathcal{P}^{m\times n}_r := \{P_0 + \lambda P_1 : P_0, P_1 \in {\mathbb C}^{m\times n}, \ \nrank(P_0 + \lambda P_1) \le r\}.$
\begin{definition} {\rm \cite[Section 2]{DeDL17} }  \label{def.homogsubsets.pen}  Let $m,n$ and $r$ be integers such that $m,n \geq 2$ and $0 < r < \min\{m,n\}$, and $a$ be an integer such that $0 \leq a \leq r$. Let us define the following subsets of $\mathcal{P}^{m\times n}_r$ 
$$
\mathcal{C}^{m\times n}_{r,a}  =  \left\{ u_1(\lambda)v_1(\lambda)^T+ \cdots +u_r(\lambda)v_r(\lambda)^T \, : \,
\begin{array}{l}
u_i(\lambda) \in \mathbb C[\lambda]^{m}, \; v_j(\lambda) \in \mathbb C[\lambda]^{n},\\
\deg (u_i) \le 1, \quad \mbox{for $i=1,\ldots ,r$}, \\
\deg (v_j) \le 1, \quad \mbox{for $j=1,\ldots ,r$}, \\
\deg (u_{1}) = \dots = \deg (u_{a}) = 0,\\
\deg (v_{a+1}) = \dots = \deg (v_{r}) = 0
\end{array}
\right\}\, .
$$
\end{definition}
We recall in Theorem \ref{lem:P_as_union_of_C's} a description of $\mathcal{P}^{m\times n}_r$  as the union of the $r+1$ sets of matrix pencils from Definition \ref{def.homogsubsets.pen}. This description is crucial for our method to compute the closest matrix pencil in $\mathcal{P}^{m\times n}_r$ to a given $m \times n$ matrix pencil $P_0 + \lambda P_1$. 
\begin{theorem}{\rm (\cite[Lemma 3.1]{DeDo16} and \cite[Lemma 4]{DeDL17})}
\label{lem:P_as_union_of_C's} Let $\mathcal{P}^{m\times n}_r$ be the set of $m\times n$ matrix pencils with normal rank at most $r$, $0<r<\min\{m,n\}$, and $\mathcal{C}^{m\times n}_{r,a}$ be the subset in Definition~\ref{def.homogsubsets.pen}, for $0 \leq a \leq r$. Then
\begin{equation}\label{P_as_union_of_C's}
\mathcal{P}^{m\times n}_r = \mathcal{C}^{m\times n}_{r,0} \cup \mathcal{C}^{m\times n}_{r,1} \cup \dots \cup \mathcal{C}^{m\times n}_{r,r}.
\end{equation}
\end{theorem}
For the sake of brevity, we do not explain the nature and properties of the description of $\mathcal{P}^{m\times n}_r$ in Theorem \ref{lem:P_as_union_of_C's}, nor its relation with the generic eigenstructures of pencils with bounded rank, but we refer the curious reader to \cite{DeDL17,DmDo17,DmDV23} for more details. 


\subsection{Formulation of the minimization problem for matrix pencils}
\label{sec:form_min_pen}
Let $P_0 + \lambda P_1$ be a given $m \times n$ complex matrix pencil. 
Our goal is to find a matrix pencil in ${\cal P}_{r}^{m \times n}$ that is nearest to $P_0 + \lambda P_1$ and the distance to it, i.e., to find
\begin{equation} \label{pro:minp}
\begin{aligned}
\dist_{{\cal P}_{r}^{m \times n}}(P_0 + \lambda P_1) &= \min_{S_0 + \lambda S_1 \in {\cal P}_{r}^{m \times n}} || (P_0-S_0) + \lambda (P_1-S_1) ||_F,\\
L_0 + \lambda L_1 &\in \argmin_{S_0 + \lambda S_1 \in {\cal P}_{r}^{m \times n}} || (P_0-S_0) + \lambda (P_1-S_1) ||_F. 
\end{aligned}
\end{equation}
The method we develop gives us a minimizer $L_0 + \lambda L_1$ and thus also the distance $\dist_{{\cal P}_{r}^{m \times n}}(P_0 + \lambda P_1)$. For brevity and convenience, in the rest of the section we will not repeat the transformations in the expression for the minimizer in \eqref{pro:minp}. 

Using Theorem \ref{lem:P_as_union_of_C's}, we can substitute the minimization problem \eqref{pro:minp} with $r+1$ minimization problems as follows 
\begin{equation*} 
\begin{aligned}
\dist_{{\cal P}_{r}^{m \times n}} \left(P_0 + \lambda P_1 \right) 
&= 
\min_{0\leq a \leq r} \left( \min_{S_0 + \lambda S_1 \in {\cal C}_{r,a}^{m\times n}} 
\left\| (P_0-S_0) + \lambda (P_1-S_1) \right\|_F \right) \\ & = 
\min_{0\leq a \leq r} \dist_{{\cal C}_{r,a}^{m\times n}}\left(P_0 + \lambda P_1 \right). \\
\end{aligned}
\end{equation*} 
Now the explicit description of a pencil $S_0 + \lambda S_1 \in {\cal C}_{r,a}^{m\times n}$, see Definition~\ref{def.homogsubsets.pen}, i.e., 
\begin{equation*}
\begin{aligned}
S_0 + \lambda S_1 = & \;
{u}_{0, 1}^0 
({v}_{0,1}^{1} + \lambda  {v}_{1, 1}^1)^{T} + \dots + 
{u}_{0,a}^0  
({v}_{0,a}^{1} + \lambda  {v}_{1, a}^1)^{T} \\
&+  ({u}_{0,a+1}^1 + \lambda  {u}_{1, a+1}^1) 
{v}_{0,a+1}^{0,T} + \dots +  
({u}_{0,r}^1 + \lambda  {u}_{1, r}^1)
{v}_{0,r}^{0,T},
\end{aligned}
\end{equation*}
in terms of arbitrary vectors $u_{0,j}^0$, $1\leq j \leq a$, $u_{0,j}^1, u_{1,j}^1$, $(a+1) \leq j \leq r$, in $\mathbb C^m$ and $v_{0,j}^0$, $(a+1) \leq j \leq r$, $v_{0,j}^1,  v_{1,j}^1$, $1\leq j \leq a$, in $\mathbb C^n$,  allows us to formulate the minimization problem for (the square of) $\dist_{{\cal C}_{r,a}^{m\times n}}\left(P_0 + \lambda P_1 \right)$, for each $a$, as the following {\em unconstrained minimization problem}
\begin{equation}
\begin{aligned}
\label{eq:totmin}
\min_{\begin{array}{c}  \scriptstyle  {u_{0,j}^0}, {u}_{0,j}^1, {u}_{1,j}^1 \in \mathbb C^m\\ \scriptstyle  {v_{0,j}^0}, {v}_{0,j}^1,  {v}_{1,j}^1 \in \mathbb C^n\end{array}}  &\left\Vert P_0 - \left( \sum_{j = 1}^a  {u}_{0,j}^0 {v}_{0,j}^{1,T} + \sum_{j = a + 1}^{r}  {u}_{0,j}^1  {v}_{0,j}^{0,T}  \right) \right\Vert_F^2 
\\  &
 + \left\Vert P_1 - \left( \sum_{j = 1}^a  {u}_{0,j}^{0} {v}_{1,j}^{1,T} + \sum_{j = a + 1}^{r}  {u}_{1,j}^{1}  {v}_{0,j}^{0,T}  \right) \right\Vert_F^2. 
\end{aligned}
\end{equation}
Equivalently, using matrix notation, we have for each $a$
\begin{equation}
\begin{aligned}
\label{eq:mainminpen}
 &(\dist_{{\cal C}_{r,a}^{m\times n}}\left(P_0 + \lambda P_1 \right))^2 \\   & \qquad = \min_{\begin{array}{c}  \scriptstyle  {U}_0^0, {U}_0^1, {U}_1^1,\\ \scriptstyle  {V}_0^0, {V}_0^1,  {V}_1^1\end{array}}  \left\|P_0-{U}_0^1 {V_0^0}^{T}-{U_0^0} {V^1_0}^T\right\|_F^2+\left\|P_1-{U}_1^1 {V_0^0}^{T}-{U_0^0} {V_1^1}^T\right\|_F^2, 
\end{aligned}
\end{equation}
where
\begin{alignat*}{3}
& {U_0^0}  = \left[ {u}_{0,1}^0 \, \cdots \, {u}_{0,a}^0 \right] \in \mathbb{C}^{m\times a}, \quad
&&{V_0^0}=\left[{v}_{0,a+1}^0 \, \cdots \, {v}_{0,r}^0\right] \in \mathbb{C}^{n\times (r-a)}, \\
& {U}_0^1 = \left[ {u}_{0,a+1}^1 \, \cdots \, {u}_{0,r}^1 \right] \in \mathbb{C}^{m\times (r-a)}, \quad
&&{V}_0^1=\left[{v}_{0,1}^1 \, \cdots \, {v}_{0,a}^1\right] \in \mathbb{C}^{n\times a}, \\
& {U}_1^1 = \left[ {u}_{1,a+1}^1 \, \cdots \, {u}_{1,r}^1 \right] \in \mathbb{C}^{m\times (r-a)}, \quad 
&& {V}_1^1=\left[{v}_{1,1}^1 \, \cdots \, {v}_{1,a}^1\right]\in \mathbb{C}^{n\times a}.
\end{alignat*}
Note that the sizes of the matrices in \eqref{eq:mainminpen} depend on $a$. In particular, if $a=0$, then the summands $U_0^0{V_0^1}^{T}$ and $U_0^0 {V_1^1}^{T}$ are absent, and if $a=r$, then the summands $U_0^1{V_0^0}^{T}$ and $U_1^1 {V_0^0}^{T}$ are absent. 

\subsection{Solution of the minimization problem over each ${\cal C}_{r,a}^{m\times n}$}
\label{sec:sol_min_pen}
An optimization problem consisting of several independent sets of unknowns may in some cases be efficiently solved in a sequence set by set, i.e., we alternate between solving a problem in which the unknowns of one set are fixed and solving a problem in which the unknowns of another set are fixed. 
The problem in (\ref{eq:mainminpen}) admits several variants of independent sets from which an alternating optimization algorithm can be formulated. 

In this work we will use the alternating least squares method, see, e.g.,~\cite{LiLR16, SrJa03}. We first summarize the idea of the method in the case of a very simple problem and then describe how it applies for the problem in \eqref{eq:mainminpen}. 
Let $R \in \mathbb{C}^{m\times n}$ be a given matrix, $S \in \mathbb{C}^{m\times p}$, $T \in \mathbb{C}^{p\times n}$ and consider the problem 
\begin{equation}
\label{eq:exmin}
\min_{S,T} \left\Vert R - ST  \right\Vert_F^2 \, .
\end{equation}
First we fix the matrix $S$ and compute the Matlab's backslash solution $\widehat{T}$ to $\displaystyle \min_{T} \left\Vert R - ST  \right\Vert_F^2$. Now we take \ $S$ \ as the unknown in \ \eqref{eq:exmin} \ and fix \ $\widehat{T}$ \ to solve the transposed problem $\displaystyle \min_{S} \left\Vert R^T -  \widehat{T}^T  S^T \right\Vert_F^2,$ obtaining the Matlab's backslash solution $\widehat{S}$. Then $\widehat{S}$ is fixed to compute with Matlab's backslash the next iteration of $T$, and so on until convergence. 

  
For our problem \eqref{eq:mainminpen}, we choose to alternate between the $U$- and $V$-variables, as this allows a nice reformulation of the problem.
We (randomly) choose the matrices $U_0^0,{U}_0^1,$ and ${U}_1^1$, resulting in the minimization over ${V}$-variables
$$
\min_{\begin{array}{c} \scriptstyle  {V}_0^0, {V}_0^1,  {V}_1^1\end{array}}  \left\|P_0-{U}_0^1 {V_0^0}^{T}-{U_0^0} {V^1_0}^T\right\|_F^2+\left\|P_1-{U}_1^1 {V_0^0}^{T}-{U_0^0} {V_1^1}^T\right\|_F^2, 
$$
which can be equivalently written as  
\begin{equation}
\label{eq:min_v_pen}
\min_{{V_0^0},{V_0^1},{V}_1^1}
\left\| 
\left[\begin{array}{l}
P_0 \\
P_1
\end{array}\right]-
\left[\begin{array}{lll}
{U}_0^1 & {U_0^0} & 0 \\
{U}_1^1 & 0 & {U_0^0}
\end{array}\right]\left[\begin{array}{l}
{V_0^0}^T \\
{V_0^1}^{T} \\
{V_1^1}^{T}
\end{array}\right]\right\|_F^2.
\end{equation}
Analogously, the minimization problem over the ${U}$-variables can be written as 
\begin{equation}
\label{eq:min_u_pen}
\min_{{U_0^0},{U_0^1},{U}_1^1} 
\left\| 
\left[\begin{array}{l}
P_0^{T} \\
P_1^{T}
\end{array}\right]-
\left[\begin{array}{lll}
{V}_0^1 & {V_0^0} & 0 \\
{V}_1^1 & 0 & {V_0^0}
\end{array}\right]\left[\begin{array}{l}
{U_0^0}^T \\
{U_0^1}^{T} \\
{U_1^1}^{T}
\end{array}\right]\right\|_F^2.
\end{equation}
So using a random starting guess for the $U$-variables and Matlab's backslash command, we solve the problem for the $V$-variables, i.e., \eqref{eq:min_v_pen}, then use this solution and Matlab's backslash command for solving the problem for the $U$-variables, i.e., \eqref{eq:min_u_pen}, and we alternate between the solution of the least squares problems \eqref{eq:min_v_pen} and \eqref{eq:min_u_pen} until convergence. See Section \ref{sec:conv} for more details on convergence. 
\begin{remark}
As mentioned at the end of Section \ref{sec:form_min_pen}, if $a=0$, then the summands $U_0^0{V_0^1}^{T}$ and $U_0^0 {V_1^1}^{T}$ disappear from \eqref{eq:mainminpen}, which in turn results in least squares problems \eqref{eq:min_v_pen} and \eqref{eq:min_u_pen} with the simpler shapes 
$$
\min_{{V_0^0}}
\left\| 
\left[\begin{array}{l}
P_0 \\
P_1
\end{array}\right]-
\left[\begin{array}{lll}
{U}_0^1 \\
{U}_1^1 
\end{array}\right]
{V_0^0}^T 
\right\|_F^2
\text{ and }
\min_{{U_0^1},{U}_1^1} 
\left\| 
\left[\begin{array}{l}
P_0^{T} \\
P_1^{T}
\end{array}\right]-
\left[\begin{array}{lll}
 {V_0^0} & 0 \\
 0 & {V_0^0}
\end{array}\right]\left[\begin{array}{l}
{U_0^1}^{T} \\
{U_1^1}^{T}
\end{array}\right]\right\|_F^2.
$$
Similarly, if $a=r$, then the summands $U_0^1{V_0^0}^{T}$ and $U_1^1 {V_0^0}^{T}$ disappear from \eqref{eq:mainminpen}, resulting in the following simplified versions of the problems  \eqref{eq:min_v_pen} and \eqref{eq:min_u_pen} 
$$
\min_{{V_0^1},{V}_1^1}
\left\| 
\left[\begin{array}{l}
P_0 \\
P_1
\end{array}\right]-
\left[\begin{array}{lll}
 {U_0^0} & 0 \\
 0 & {U_0^0}
\end{array}\right]\left[\begin{array}{l}
{V_0^1}^{T} \\
{V_1^1}^{T}
\end{array}\right]\right\|_F^2
\text{ and }
\min_{{U_0^0}} 
\left\| 
\left[\begin{array}{l}
P_0^{T} \\
P_1^{T}
\end{array}\right]-
\left[\begin{array}{lll}
{V}_0^1 \\
{V}_1^1 
\end{array}\right]
{U_0^0}^T 
\right\|_F^2.
$$
\end{remark}

\begin{remark}  
Observe that
\begin{equation} \label{eq.sizesuvmats}
\left[\begin{array}{lll}
{U}_0^1 & {U_0^0} & 0 \\
{U}_1^1 & 0 & {U_0^0}
\end{array}\right] \in \mathbb{C}^{2m \times (r+a)} \; \text{and} \;
\left[\begin{array}{lll}
{V}_0^1 & {V_0^0} & 0 \\
{V}_1^1 & 0 & {V_0^0}
\end{array}\right] \in \mathbb{C}^{2n \times (2r-a)}.
\end{equation}
So, both matrices have more rows than columns because $0\leq a \leq r < \min\{m,n \}$. Therefore, generically, the least squares problems \eqref{eq:min_v_pen} and \eqref{eq:min_u_pen} have a unique solution. In the non-generic case that the columns of some of the matrices in \eqref{eq.sizesuvmats} are linearly dependent, the corresponding least squares problem has infinitely many solutions. Among these infinitely many solutions we take the one given by Matlab's backslash command, which is not the one provided by the Moore-Penrose inverse, in general. We emphasize that the convergence results in Section \ref{sec:conv} are valid for any of the infinitely many possible solutions under the natural assumption that the sequence of solutions of the alternating least squares problems remains bounded.
\end{remark}

\section{Approximation with low-rank matrix polynomials} \label{sec:pol}

Similarly to Section \ref{sec:pencase} for matrix pencils, the main goal of this section is to formulate the minimization problem for matrix polynomials of arbitrary grade in terms of a finite number of minimization sub-problems over certain sets of matrix polynomials that can be explicitly described in factorized form and to solve each of these subproblems using the alternating least squares method. The main theoretical ingredients are now Definition \ref{def.homogsubsets} and Theorem~\ref{thm.relationgeneigen}, which are the counterparts of Definition \ref{def.homogsubsets.pen} and Theorem \ref{lem:P_as_union_of_C's} used in the case of pencils.

\subsection{Preliminary information about matrix polynomials}

Recall that a matrix polynomial $P(\lambda ) = \sum_{i= 0}^d \lambda^iP_i$ has a grade $d$, i.e., a degree less than or equal to $d$. If $P_d \neq 0$, then $P(\lambda)$ is a matrix polynomial with degree $d$. 
In addition to the notation introduced at the beginning of Subsection \ref{sec.prelimpenc}, $\mathbb{C}[\lambda]^{m \times n}_{d}$ denotes the set of $m\times n$ matrix polynomials with grade $d$ and $\mathbb{C}[\lambda]^{m \times n}_{d,r}$ denotes the set of $m\times n$ matrix polynomials with grade $d$ and normal rank at most $r$. Notably, $\mathbb{C}[\lambda]^{m \times n}_{d} \supseteq \mathbb{C}[\lambda]^{m \times n}_{d,r}$, where equality holds if $r = \min\{m,n\}$. 
A distance on $\mathbb{C}[\lambda]^{m \times n}_{d}$ is defined as $\rho(P,Q) = \|P(\lambda) - Q(\lambda)\|_F$, making $\mathbb{C}[\lambda]^{m \times n}_{d}$ a metric space and allowing us to define the closure of any subset $\mathcal{S} \subseteq \mathbb{C}[\lambda]^{m \times n}_{d}$, which is denoted as $\overline{\mathcal{S}}$.
Given a matrix polynomial $P(\lambda)$, its $i$th row is represented by $P_{i*} (\lambda)$, or simply by $P_{i*}$, and its $j$th column by $P_{*j} (\lambda)$, or simply by $P_{*j}$. 
The following subsets of $\mathbb{C}[\lambda]^{m \times n}_{d}$ are fundamental in this work.
\begin{definition} {\rm \cite[Definition 4.12]{DmDV23} } \label{def.homogsubsets} Let $m,n,r$ and $d$ be integers such that $m,n \geq 2$, $d \geq 1$ and $0 < r < \min\{m,n\}$, and $a$ be an integer such that $0 \leq a \leq rd$. Let us define $d_R := \left\lfloor a /r \right\rfloor$, $t_R := a \; \mathrm{mod} \; r$ and the following subset of polynomial matrices 
$$
\mathcal{C}^{m\times n}_{d,r,a}  :=  \left\{ L(\lambda) R(\lambda) \, : \,
\begin{array}{l}
L(\lambda) \in \mathbb C[\lambda]^{m\times r}, \; R(\lambda) \in \mathbb C[\lambda]^{r\times n},\\
\deg (R_{i*}) \leq d_R + 1, \quad \mbox{for $i=1,\ldots ,t_R$}, \\
\deg (R_{i*}) \leq d_R, \quad \mbox{for $i=t_R + 1,\ldots ,r$}, \\
\deg (L_{*i}) \leq d - d_R - 1, \quad \mbox{for $i=1,\ldots ,t_R$}, \\
\deg (L_{*i}) \leq d - d_R, \quad \mbox{for $i= t_R + 1,\ldots ,r$}
\end{array}
\right\}  \subset \mathbb C[\lambda]^{m\times n}_{d,r} \, .
$$
\end{definition}
Note that for $d=1$, we have that $\mathcal{C}^{m\times n}_{1,r,a} = \mathcal{C}^{m\times n}_{r,a}$, with $\mathcal{C}^{m\times n}_{r,a}$ the set of pencils in Definition \ref{def.homogsubsets.pen}. 

The following theorem portrays $\mathbb{C}[\lambda]^{m \times n}_{d,r}$, for $0< r < \min\{m,n\}$ and $d\ge 1$, as the union of the closures of the $rd+1$ sets of matrix polynomials from Definition \ref{def.homogsubsets}.

\begin{theorem}\label{thm.relationgeneigen} {\rm \cite[Theorem 4.13-(iii)]{DmDV23} } Let $\mathcal{C}^{m\times n}_{d,r,a}$ be the set of matrix polynomials introduced in Definition \ref{def.homogsubsets}. Then
 $$\displaystyle \mathbb C[\lambda]^{m\times n}_{d,r} = \bigcup_{0 \leq a \leq rd} \overline{\mathcal{C}^{m\times n}_{d,r,a}}.$$
\end{theorem} 
Observe that Theorem \ref{lem:P_as_union_of_C's} is stronger than the result obtained from making $d=1$ in Theorem \ref{thm.relationgeneigen}, because Theorem \ref{lem:P_as_union_of_C's} does not involve closures in the sets $\mathcal{C}^{m\times n}_{r,a}$. A detailed discussion about this issue can be found in \cite[Remark 4.14-Example 4.15]{DmDV23}. Once again we refer the reader to \cite{DmDo17,DmDV23} for more insight on the description of $\displaystyle \mathbb C[\lambda]^{m\times n}_{d,r}$ in Theorem \ref{thm.relationgeneigen}. 

\subsection{Formulation of the minimization problem for matrix polynomials}
%

Let $P(\lambda) = \sum_{i= 0}^d \lambda^iP_i \in \mathbb C[\lambda]^{m\times n}_{d}$.  Our  goal is to find a matrix polynomial in $\mathbb C[\lambda]^{m\times n}_{d,r}$ that is nearest to $P(\lambda)$ and the distance to it, i.e.,  
\begin{equation} \label{pro:polr_gen}
\begin{aligned}
\dist_{\mathbb C[\lambda]^{m\times n}_{d,r}} \left(\sum_{i= 0}^d \lambda^iP_i \right) &= \min_{\sum_{i= 0}^d \lambda^iS_i \in \mathbb C[\lambda]^{m\times n}_{d,r}} \left\| \sum_{i= 0}^d \lambda^i(P_i-S_i) \right\|_F, \\
\sum_{i= 0}^d \lambda^iL_i &\in \argmin_{\sum_{i= 0}^d \lambda^iS_i \in \mathbb C[\lambda]^{m\times n}_{d,r}} \left\|  \sum_{i= 0}^d \lambda^i(P_i-S_i) \right\|_F .
\end{aligned}
\end{equation}
As in the case of pencils, we will not repeat the transformations in the expression for the minimizer in \eqref{pro:polr_gen} but we emphasize that the method we develop gives us a minimizer $\sum_{i= 0}^d \lambda^iL_i$ and thus also the distance $\dist_{\mathbb C[\lambda]^{m\times n}_{d,r}} \left(\sum_{i= 0}^d \lambda^iP_i \right)$. 

Similarly as we did for matrix pencils in Section \ref{sec:form_min_pen}, but now using Theorem \ref{thm.relationgeneigen}, we substitute the minimization problem \eqref{pro:polr_gen} with $rd+1$ minimization problems as follows
\begin{equation} \label{pro:polr}
\begin{aligned}
\dist_{\mathbb C[\lambda]^{m\times n}_{d,r}} \left(\sum_{i= 0}^d \lambda^iP_i \right) 
& = \min_{ \sum\limits_{i= 0}^d \lambda^iS_i \in \left(\bigcup\limits_{0 \leq a \leq rd} \overline{\mathcal{C}^{m\times n}_{d,r,a}}\right)} \left\| \sum_{i= 0}^d \lambda^i(P_i-S_i) \right\|_F \\ 
&= 
\min_{0\leq a \leq rd} \left( \min_{\sum\limits_{i= 0}^d \lambda^iS_i \in \mathcal{C}^{m\times n}_{d,r,a}} \left\| \sum_{i= 0}^d \lambda^i(P_i-S_i) \right\|_F \right) \\ & = \min_{0\leq a \leq rd} \dist_{\mathcal{C}^{m\times n}_{d,r,a}}\left(\sum_{i= 0}^d \lambda^iP_i \right).
\end{aligned}
\end{equation} 
 Now we consider each of the $rd+1$ minimization problems derived above: 
\begin{equation} \label{pro:minc}
\begin{aligned}
&\dist_{\mathcal{C}^{m\times n}_{d,r,a}}\left(\sum_{i= 0}^d \lambda^iP_i \right) 
= \min_{\sum\limits_{i= 0}^d \lambda^iS_i \in {\mathcal{C}^{m\times n}_{d,r,a}}} \left\| \left[ \begin{array}{c}
P_0 \\
P_1  \\
 \vdots \\
P_{d}  \\
\end{array}\right] - \left[ \begin{array}{c}
S_0 \\
S_1  \\
 \vdots \\
S_{d}  \\
\end{array}\right] \right\|_F.
\end{aligned}
\end{equation} 
The remainder of this subsection is devoted to reformulating \eqref{pro:minc} in a form that explicitly exploits the structure of the polynomials in $\mathcal{C}^{m\times n}_{d,r,a}$; see Definition~\ref{def.homogsubsets}. For this, we introduce some additional definitions and notation. Since we consider $rd+1$ cases, one for each value of $a$, our goal is to develop a unified approach that covers all of them. Example~\ref{ex:deg2} is included to illustrate these definitions in a simple setting.  

For given integers $d,m,n,r,a$, with $d \geq 1, m,n \geq 2, 0 < r < \min \{m,n\}, 0 \leq a \leq rd$, we define two partitioned matrices $U_a$ and $V_a$ as follows
\begin{align} \label{eq.uacol}
U_a^T &= [ U_0^{0} \ | \ U_0^{1} \ U_1^{1} \ | \ U_0^{2} \ U_1^{2} \ U_2^{2} \ | \cdots | \ U_0^{d-1} \cdots U_{d-1}^{d-1} | \ U_0^{d} \cdots U_d^{d} ], \quad \text{ and } \\ \label{eq.vacol}
V_{a}^T &= [ V_0^{0} \ | \ V_0^{1} \ V_1^{1} \ | \ V_0^{2} \ V_1^{2} \ V_2^{2} \ | \cdots | \ V_0^{d-1} \cdots V_{d-1}^{d-1} | \ V_0^{d} \cdots V_d^{d} ],
\end{align}
where $U_{i}^{d-k} \in \mathbb C^{m\times n_{d-k}},  V_{j}^{k} \in \mathbb C^{n\times n_{d-k}}, n_{d-k}:= \max \{0, r-|a-kr| \}, \ k = 0, 1, \dots , d $, $ i=0, 1, \dots , d-k$, and $ j=0, 1, \dots , k$. Note that $n_{d-k}$ is the number of columns in each of the matrices $U_{i}^{d-k}$ and $V_{j}^{k}$ and that if $n_{d-k}=0$, then the corresponding matrices are absent from $U_a$ and $V_a$. Next, we define the matrix
\begin{equation} \label{eq.matua}
    \begin{aligned}
&{\bf M}(U_a) := \begin{bmatrix}
    M_0 (U_a) | M_1 (U_a) | M_2 (U_a) | \cdots | M_{d-1} (U_a)| M_{d} (U_a)
\end{bmatrix} \\ 
& :=
\scriptsize
\left[ 
\setlength{\arraycolsep}{3.5pt}
\scalemath{1.0}{
\begin{array}{c|c c| c c c | c | c c c c | c c c c c c}
U_0^d & U_0^{d-1} & 0 & U_0^{d-2} & 0 & 0 & & U_0^{1} & & 0 & 0 &  U_0^{0} &&&&& \\
U_1^d & U_1^{d-1} & U_0^{d-1} & \vdots & U_0^{d-2} & 0 & & U_1^{1} &  & \vdots & 0 & & U_0^{0} &&&& \\
\vdots & \vdots & U_1^{d-1} & U_{d-3}^{d-2} & \vdots & U_0^{d-2} & \dots & 0 &\ddots & 0 & \vdots & && U_0^{0} &&& \\
\vdots & U_{d-2}^{d-1} &\vdots &  U_{d-2}^{d-2} & U_{d-3}^{d-2} & \vdots & & 0 &  & U_0^{1} & 0 & &&& \ddots && \\
U_{d-1}^d & U_{d-1}^{d-1} &U_{d-2}^{d-1} &  0 & U_{d-2}^{d-2} & U_{d-3}^{d-2} & & \vdots & & U_1^{1} & U_0^{1} & &&&& U_0^{0} & \\
U_d^d & 0 &U_{d-1}^{d-1} &  0 & 0 & U_{d-2}^{d-2} & & 0 & & 0 & U_1^{1} & &&&&& U_0^{0} \\
\end{array} }\right] .
    \end{aligned}
\end{equation}
Observe that $M_k (U_a) \in \mathbb{C}^{(d+1)m \times (k+1) n_{d-k}}$ has the same number of columns as $[ V_0^{k} \cdots V_{k}^{k} ] \in \mathbb{C}^{n \times (k+1) n_{d-k}}$ for $k=0, 1,\ldots,d$. Analogously, we define 
\begin{align}
{\bf M}(V_a) & := \begin{bmatrix}
    M_0 (V_a) | M_1 (V_a) | M_2 (V_a) | \cdots | M_{d-1} (V_a)| M_{d} (V_a)
\end{bmatrix} \label{eq.matva},
\end{align}
by replacing in \eqref{eq.matua} each matrix $U^j_i$ by $V^j_i$. In this case, the sizes of the submatrices are as follows: $M_k (V_a) \in \mathbb{C}^{(d+1)n \times (k+1) n_{k}}$ and $[ U_0^{k} \cdots U_{k}^{k} ]  \in \mathbb{C}^{m \times (k+1) n_{k}}$ for $k=0, 1,\ldots,d$.



Lemma \ref{ndklemma} shows that $n_{d-k}= \max \{0, r-|a-kr| \}$ is strictly greater than zero for at most two values of $k$. Thus, among $U_0^{0},U_0^{1},U_0^{2}, \dots , U_0^{d}$ at most two $U_{0}^{d-k}$ are present, i.e., they have a positive number of columns. Clearly, the absence of $U_{0}^{d-k}$ implies the absence of $U_{i}^{d-k}$, for $i=0, \dots , d-k$. Similarly, among $V_0^{0},V_0^{1},V_0^{2},\dots , V_0^{d}$ in \eqref{eq.vacol} at most two $V_{0}^{k}$ are present.
\begin{lemma}
\label{ndklemma}
Let $d,r,$ and $a$ be given integers, with $d$ and $r$ positive and such that $0 \leq a \leq rd$. Consider $n_{d-k}= \max \{0, r-|a-kr| \}$ for $ k = 0, 1, \dots, d$. Then, the following two statements hold:  
\begin{itemize}
\item[\rm (a)] $n_{d-k} > 0$ for at most two values of $k$, namely for $k \in (a/r -1, a/r +1) \cap \mathbb Z$;
\item[\rm (b)] $\sum_{k=0}^d n_{d-k} =r$. 
\end{itemize}   
\end{lemma}
\begin{proof}
    The solution of $r-|a-kr|>0$ shows the first item. For the second item, note that $k \in (a/r -1, a/r +1) \cap \mathbb Z$ means $k = a/r$ if $a/r = \lfloor a/r \rfloor$, and $k$ is equal to $\lfloor a/r \rfloor$ and to $\lfloor a/r \rfloor+1$ if $a/r \neq \lfloor a/r \rfloor$. Adding $n_{d-k}$ for these values of $k$ proves the second item. 
\end{proof}

As a corollary of Lemma \ref{ndklemma}, its proof and the number of columns of the submatrices of \eqref{eq.uacol}-\eqref{eq.matva}, we obtain the following result about the matrices ${\bf M} (U_a)$,  $V_a$, ${\bf M} (V_a)$, and  $U_a$. Its simple proof is omitted.
\begin{lemma} \label{lemm.sizesuava} Let $m, n, r, d$ and $a$ be integers such that $m,n \geq 2$, $d\geq 1$, $0 < r < \min \{m,n\}$ and $0 \leq a \leq rd$. Let $d_R = \lfloor a/r \rfloor$, $t_R =  a \; \mathrm{mod} \; r$, $n_{d-k}= \max \{0, r-|a-kr| \}$,  ${\bf M} (U_a)$ and $V_a$ be the matrices in \eqref{eq.matua} and \eqref{eq.vacol}, and ${\bf M} (V_a)$ and $U_a$ be the matrices in \eqref{eq.matva} and \eqref{eq.uacol}. 
\begin{itemize}
    \item[\rm (a)] If $d_R = a/r$, then $n_{d-d_R} = r$, 
    \begin{align*}
    & {\bf M}(U_a) = M_{d_R} (U_a) \in \mathbb{C}^{(d+1)m \times (r+a)}, \\
    & V_a^T = [V_0^{d_R} \ V_1^{d_R} \cdots V_{d_R}^{d_R} ] \in \mathbb{C}^{n \times (r+a)}, \\
    & {\bf M}(V_a) = M_{d-d_R} (V_a) \in \mathbb{C}^{(d+1)n \times ((d+1)r-a)}, \quad \text{and} \\
    & U_a^T = [U_0^{d-d_R} \ U_1^{d-d_R} \cdots U_{d-d_R}^{d-d_R} ]  \in \mathbb{C}^{m \times ((d+1)r-a)}. 
    \end{align*}

    \item[\rm (b)] If $d_R \ne a/r$, then $n_{d-d_R} = r-t_R$, $n_{d-d_R-1} = t_R$, 
    \begin{align*}
    & {\bf M}(U_a) =  [M_{d_R} (U_a) | M_{d_R+1} (U_a)] \in \mathbb{C}^{(d+1)m \times (r+a)}, \\  
    &V_a^T =  [ V_0^{d_R}\ V_1^{d_R} \cdots V_{d_R}^{d_R} | \ V_0^{d_R +1} \  V_1^{d_R +1}  \cdots V_{d_R +1}^{d_R + 1} ] \in \mathbb{C}^{n \times (r+a)}, \\
    & {\bf M}(V_a) =  [M_{d-d_R-1} (V_a) | M_{d-d_R} (V_a)] \in \mathbb{C}^{(d+1)n \times ((d+1)r-a)}, \quad \text{and} \\  
    &U_a^T =  [ U_0^{d-d_R-1} \ U_1^{d-d_R-1} \cdots U_{d-d_R-1}^{d-d_R-1} | \ U_0^{d - d_R} \ U_1^{d - d_R} \cdots U_{d-d_R}^{d-d_R} ] \in \mathbb{C}^{m \times ((d+1)r-a)}. \\
    \end{align*}
    
\end{itemize} 
\end{lemma}

\begin{example}   
\label{ex:deg2}
We illustrate the definitions of ${\bf M} (U_a)$ and $V_a$ in the case $d=2$ and $r=2$ for all the possible values of $a$. 
In this situation
\begin{align} \label{eq.uacol_2}
U_a^T &= [U_0^{0} \ | \ U_0^{1}  \ U_1^{1}  \ | \ U_0^{2}  \   U_1^{2}  \  U_2^{2}], \ \ 0 \leq a \leq 4,
\end{align}
before removing the submatrices with zero columns. The corresponding full block matrix reads
\begin{equation*}
{\bf M}(U_a) = \begin{bmatrix}
    M_0 (U_a) | M_1 (U_a) | M_{2} (U_a)
\end{bmatrix} =
\left[ \begin{array}{c|c c| c c c }
U_0^2 & U_0^{1} & 0 & U_0^{0} & 0 & 0  \\
U_1^2 & U_1^{1} & U_0^{1} & 0 & U_0^{0} & 0  \\
U_2^2 & 0 & U_1^{1} & 0 & 0 & U_0^{0}  \\
\end{array}\right]
\end{equation*}
and
\begin{equation*}
    V_a^T = [ V_0^{0} \  | \ V_0^{1} \ V_1^{1} \ | \ V_0^{2} \ V_1^{2} \ V_2^{2} ]. 
\end{equation*}
For each $a$, the matrices ${\bf M}(U_a)$ and $V_a^T$ and their sizes are obtained from Lemma \ref{lemm.sizesuava}. The results are summarized in Table~\ref{tab:deg2_cases}. The matrices ${\bf M} (V_a)$ and $U_a^T$ can be obtained analogously for $a= 0,1,2,3,4$. For brevity, we just show the case $a=3$:
$$ {\bf M} (V_3) =
\left[
\begin{array}{c|cc}
V_0^2 & V_0^{1} & 0 \\
V_1^2 & V_1^{1} & V_0^{1} \\
V_2^2 & 0 & V_1^{1}
\end{array}
\right] \in \mathbb{C}^{3n \times 3}, \quad  U_3^T = [\, U_0^{0} \mid U_0^{1}\ U_1^{1} \,] \in \mathbb{C}^{m \times 3}.
$$

\begin{table}[h]
\centering
\renewcommand{\arraystretch}{0.99}
\begin{tabular}{c|l|c|l|c}
$a$ & ${\bf M}(U_a)$& size(${\bf M}(U_a)$) & $V_a^T$ & size($V_a^T$)\\
\hline \rule{0pt}{6ex}

0 &
$\displaystyle
\left[
\begin{array}{c}
U_0^2 \\
U_1^2 \\
U_2^2
\end{array}
\right]
$
&
$3m \times 2$
&
$[\, V_0^{0} \,]$
& 
$n \times 2$
\\[5ex]

1 &
$\displaystyle
\left[
\begin{array}{c|cc}
U_0^2 & U_0^{1} & 0 \\
U_1^2 & U_1^{1} & U_0^{1} \\
U_2^2 & 0 & U_1^{1}
\end{array}
\right]
$
&
$3m \times 3$
&
$[\, V_0^{0} \mid V_0^{1}\ V_1^{1} \,]$
&
$n \times 3$
\\[5ex]

2 &
$\displaystyle
\left[
\begin{array}{cc}
U_0^{1} & 0 \\
U_1^{1} & U_0^{1} \\
0 & U_1^{1}
\end{array}
\right]
$
&
$3m \times 4$
&
$[\, V_0^{1}\ V_1^{1} \,]$
&
$n \times 4$
\\[5ex]

3 &
$\displaystyle
\left[
\setlength{\arraycolsep}{3.3pt}
\scalemath{1.0}{
\begin{array}{cc|ccc}
U_0^{1} & 0 & U_0^{0} & 0 & 0 \\
U_1^{1} & U_0^{1} & 0 & U_0^{0} & 0 \\
0 & U_1^{1} & 0 & 0 & U_0^{0}
\end{array}}
\right]
$
&
$3m \times 5$
&
$
\setlength{\arraycolsep}{3.3pt}
\scalemath{1.0}{
[\, V_0^{1} V_1^{1} \mid V_0^{2} V_1^{2} V_2^{2} \,]
}$
&
$n \times 5$
\\[5ex]

4 &
$\displaystyle
\left[
\begin{array}{ccc}
U_0^{0} & 0 & 0 \\
0 & U_0^{0} & 0 \\
0 & 0 & U_0^{0}
\end{array}
\right]
$
&
$3m \times 6$
&
$[\, V_0^{2}\ V_1^{2}\ V_2^{2} \,]$
&
$n \times 6$

\end{tabular}
\smallskip
\caption{Explicit shapes and sizes of ${\bf M}(U_a)$ and $V_a^T$ for all possible $a$, in the case of $d=2$ and $r=2$.}
\label{tab:deg2_cases}
\end{table}
\qed
\end{example}

\medskip

\begin{remark} \label{rem.ourspaces}
Lemma \ref{lemm.sizesuava} shows that, for each $a$, $U_a \mapsto {\bM} (U_a)$ is a linear function of $U_a$ with domain the vector space $\mathbb{C}^{((d+1)r -a) \times m}$ and codomain  $\mathbb{C}^{(d+1) m \times (r+a)}$. Analogously, $V_a \mapsto {\bM} (V_a)$ is a linear function of $V_a$ with domain the vector space $\mathbb{C}^{(r+a) \times n}$ and codomain  $\mathbb{C}^{(d+1) n \times ((d+1)r -a)}$. This fact will be used in the sequel and, in particular, in Section \ref{sec:conv}.
\end{remark}

To connect ${\bf M}(U_a)$ and $V_a$  with the set ${\mathcal{C}^{m\times n}_{d,r,a}}$ in the minimization problem \eqref{pro:minc}, we define the set 
\begin{equation}
    \label{Btilde}
    \widetilde{\cC}^{m\times n}_{d,r,a} := \left\{
\left[ \begin{array}{c}
S_0 \\
S_1  \\
 \vdots \\
S_{d}  \\
\end{array}\right] \, : \,
    \sum_{i= 0}^d \lambda^iS_i \in {\mathcal{C}^{m\times n}_{d,r,a}} 
    \right\} \subset \mathbb{C}^{(d+1)m \times n}. 
\end{equation}
Observe that $\widetilde{\cC}^{m\times n}_{d,r,a}$ is isometric to $\cC^{m\times n}_{d,r,a}$ via the Frobenius norm. Lemma \ref{lemm.polconnec} establishes the  connection and relates the matrices ${\bf M}(U_a) \cdot V_a \in \mathbb{C}^{(d+1)m \times n}$ and 
${\bf M}(V_a) \cdot U_a \in \mathbb{C}^{(d+1)n \times m}$.


\begin{lemma} \label{lemm.polconnec}
 Let $m,n,r$ and $d$ be integers such that $m,n \geq 2$, $d \geq 1$ and $0 < r < \min\{m,n\}$, and $a$ be an integer such that $0 \leq a \leq rd$.  Let ${\bf M}(U_a)$, $V_a$, ${\bf M}(V_a)$, $U_a$ be matrices as in Lemma \ref{lemm.sizesuava}, and  $\widetilde{\cC}^{m\times n}_{d,r,a}$ be the set in \eqref{Btilde}. Then 
     \begin{equation} \label{eq.polconnec1}
        \widetilde{\mathcal{C}}^{m\times n}_{d,r,a} = \left\{  {\bf M}(U_a) \cdot V_a \ : \ U_a \in \mathbb{C}^{((d+1)r -a) \times m}, \, V_a \in \mathbb{C}^{(r +a) \times n}  \right\} .
     \end{equation}
Moreover,
\begin{equation} \label{eq.reluvvu}
\left[ \begin{array}{c}
S_0 \\
S_1  \\
 \vdots \\
S_{d}  \\
\end{array}\right] = {\bf M}(U_a) \cdot V_a \;\; \text{if and only if} \;\;
\left[ \begin{array}{c}
S_0^T \\
S_1^T  \\
 \vdots \\
S_{d}^T  \\
\end{array}\right] = {\bf M}(V_a) \cdot U_a.
\end{equation}
 \end{lemma}
 \begin{proof}
     Fix $a$ and let $[S_0^T, S_1^T, \dots, S_{d}^T]^T \in \widetilde{\mathcal{C}}^{m\times n}_{d,r,a}$, which is equivalent to $\sum_{i= 0}^d \lambda^iS_i \in {\mathcal{C}^{m\times n}_{d,r,a}}$ by the definition of $\widetilde{\mathcal{C}}^{m\times n}_{d,r,a}$. Let $d_R = \lfloor a/r \rfloor$ and assume $d_R  \neq  a/r$; the case $d_R = a/r $ is analogous but simpler, and is omitted. By Definition \ref{def.homogsubsets}, $\sum_{i= 0}^d \lambda^iS_i \in {\mathcal{C}^{m\times n}_{d,r,a}}$ is equivalent to 
     \begin{equation} \label{eq.def2explicit}
     \sum_{i= 0}^d \lambda^iS_i = 
     \left(\sum_{i= 0}^{d-d_R} \lambda^iU_i^{d-d_R}\right)
     \left(\sum_{i= 0}^{d_R} \lambda^i{V_i^{d_R}}^T\right) 
      + \left(\sum_{i= 0}^{d-d_R-1} \lambda^iU_i^{d-d_R-1}\right)
     \left(\sum_{i= 0}^{d_R+1} \lambda^i{V_i^{d_R+1}}^T\right), 
    \end{equation}
    with $\left(\sum_{i= 0}^{d-d_R} \lambda^iU_i^{d-d_R}\right) \in \mathbb{C}[\lambda]^{m \times (r-t_R)}$, $\left(\sum_{i= 0}^{d-d_R-1} \lambda^iU_i^{d-d_R-1}\right)  \in \mathbb{C}[\lambda]^{m \times t_R}$,  $\left(\sum_{i= 0}^{d_R} \lambda^i{V_i^{d_R}}^T\right) \allowbreak \in \mathbb{C}[\lambda]^{(r-t_R) \times n}$ and $\left(\sum_{i= 0}^{d_R+1} \lambda^i{V_i^{d_R+1}}^T\right) \in \mathbb{C}[\lambda]^{t_R \times n}$. The equality \eqref{eq.def2explicit} can be rewritten using the block-Toeplitz matrices $M_{d_R} (U_a)$ and $M_{d_R+1} (U_a)$ built as in \eqref{eq.matua} from the matrix polynomials $\sum_{i= 0}^{d-d_R} \lambda^iU_i^{d-d_R}$ and $\sum_{i= 0}^{d-d_R-1} \lambda^iU_i^{d-d_R-1}$, respectively, as follows
    \begin{equation*}
    \begin{aligned}
\left[ \begin{array}{c}
S_0 \\
S_1  \\
 \vdots \\
S_{d}  \\
\end{array}\right] & = M_{d_R} (U_a)
  \left[ \begin{array}{c}
{V_0^{d_R}}^T \\
{V_1^{d_R}}^T  \\
 \vdots \\
{V_{d_R}^{d_R} }^T \\
\end{array}\right] + M_{d_R+1} (U_a)
  \left[ \begin{array}{c}
{V_0^{d_R+1}}^T \\
{V_1^{d_R+1}}^T  \\
 \vdots \\
{V_{d_R+1}^{d_R+1}}^T  \\
\end{array}\right] = {\bf M}(U_a) \cdot V_a,
    \end{aligned}
\end{equation*}
where the last equality follows from Lemma \ref{lemm.sizesuava}-(b). This yields \eqref{eq.polconnec1}. Taking transposes in \eqref{eq.def2explicit} and proceeding in a similar way yields \eqref{eq.reluvvu}.
\end{proof}

Lemma \ref{lemm.polconnec} together with Remark \ref{rem.ourspaces} allow us to formulate the minimization problem \eqref{pro:minc} as the following {\em unconstrained minimization problem}
\begin{equation} \label{eq.unsconpoly1}
\begin{aligned}
&\dist_{\mathcal{C}^{m\times n}_{d,r,a}}\left(\sum_{i= 0}^d \lambda^iP_i \right) 
= \min_{ \begin{array}{c}  \scriptstyle U_a \in \mathbb{C}^{((d+1)r -a) \times m},\\ \scriptstyle  V_a \in  \mathbb{C}^{(r+a) \times n} \end{array} } 
\left\| \left[ \begin{array}{c}
P_0 \\
P_1  \\
 \vdots \\
P_{d}  \\
\end{array}\right] - {\bf M}(U_a) \cdot V_a  \right\|_F . 
\end{aligned}
\end{equation} 
It is important to note that the cost function of this minimization problem admits another alternative expression, because \eqref{eq.reluvvu} implies 
\begin{equation} \label{eq.costfunctwoexp}
\left\| \left[ \begin{array}{c}
P_0 \\
P_1  \\
 \vdots \\
P_{d}  \\
\end{array}\right] - {\bf M}(U_a) \cdot V_a  \right\|_F =     
\left\| \left[ \begin{array}{c}
P_0^T \\
P_1^T \\
 \vdots \\
P_{d}^T  \\
\end{array}\right] - {\bf M}(V_a) \cdot U_a  \right\|_F .
\end{equation}
These two expressions play a key role in the solution of the problem via the method in Section \ref{sec:sol_min_pol}.



\subsection{Solution of the minimization problem over each $\mathcal{C}^{m\times n}_{d,r,a}$}
\label{sec:sol_min_pol}
Similarly as for matrix pencils, we solve \eqref{eq.unsconpoly1} by using an alternating least squares method between the $U-$ and $V-$ variables. To be exact, we (randomly) choose a matrix $U_a \in \mathbb{C}^{((d+1)r -a) \times m}$, construct with its (adequate) submatrices $\bM (U_a)$, and solve the standard least squares problem
\begin{equation}
\label{alt_pol_min1}
\min_{V_a \in \mathbb{C}^{(r+a) \times n}}
\left\| 
\left[\begin{array}{c}
P_0 \\
\vdots \\
P_d
\end{array}\right]-
{\bf M}(U_a) \cdot V_a \right\|_F,
\end{equation}
using, in practice, Matlab's backslash command. From the computed $V_a$, we construct $\bM (V_a)$, and, taking into account \eqref{eq.costfunctwoexp}, solve the standard least squares problem
\begin{equation}
\label{alt_pol_min2}
\min_{U_a  \in \mathbb{C}^{((d+1)r -a) \times m} } 
\left\| 
\left[\begin{array}{c}
P_0^{T} \\
\vdots  \\
P_d^{T}
\end{array}\right]-
{\bf M}(V_a) \cdot U_a \right\|_F,
\end{equation}
using, in practice, again Matlab's backslash command. Then we again repeat the first step with the computed $U_a$ and so on until convergence.

\begin{remark} Observe that ${\bf M}(U_a) \in \mathbb{C}^{(d+1)m \times (r+a)}$ and
${\bf M}(V_a) \in \mathbb{C}^{(d+1)n \times ((d+1)r-a)}$ are both matrices with more rows than columns because $0 \leq a \leq rd$ and $0 < r < \min \{m,n\}$. Therefore, generically, the least squares problems \eqref{alt_pol_min1} and \eqref{alt_pol_min2} have a unique solution, as in the case of pencils. In the non-generic cases in which the columns of ${\bf M}(U_a)$ or
${\bf M}(V_a)$ are linearly dependent, the corresponding least squares problem has infinitely many solutions and we take the one given by Matlab's backslash command. The convergence analysis in Section \ref{sec:conv} is valid for any other selection of a solution under the natural assumption that the sequence of solutions of the alternating least squares problems remains bounded.  \label{rem.uniquesizemat}
\end{remark}

\begin{remark} {\rm (On the real minimization problem).} If the polynomial $P(\lambda) = \sum_{i=0}^d \lambda^i P_i \in \mathbb{R} [\lambda]^{m\times n}_d$ is real and the initial iterate $U_a \in \mathbb{R}^{((d+1)r -a) \times m}$ is chosen to be real, then all the iterates computed by the alternating least squares method described above remain real. In fact, all the results presented in Section \ref{sec:pol} are valid if $\mathbb{C}$ is replaced by $\mathbb{R}$. Therefore, the presented method solves the real version of the problem \eqref{pro:polr_gen}. \label{rem.realminization}
\end{remark}

\section{Convergence analysis}
\label{sec:conv}

In this section, we analyze the convergence of the algorithm described in Section \ref{sec:sol_min_pol}. We will work with complex matrix polynomials, but we emphasize that the analysis remains valid in the real case (recall Remark \ref{rem.realminization}). 

Let us start by introducing some notation. For brevity, we define $p_a = (d+1)r-a$ and $q_a = r+a$. According to \eqref{eq.unsconpoly1}, we take the matrices $U_a \in \mathbb{C}^{p_a \times m}$ and $V_a \in \mathbb{C}^{q_a \times n}$ in Lemma \ref{lemm.sizesuava} as the variables of the minimization problem. Recall that according to Remark \ref{rem.ourspaces}, for $\alpha, \beta \in \mathbb{C}$, $U_a, E_a\in \mathbb{C}^{p_a \times m}$, and $V_a, F_a \in \mathbb{C}^{q_a \times n}$, the equalities
\begin{equation}
\bM (\alpha \, U_a + \beta \, E_a) = \alpha \, \bM (U_a) + \beta \, \bM (E_a)  \quad \mbox{and} \quad \bM (\alpha \, V_a + \beta F_a) = \alpha \, \bM (V_a) + \beta \, \bM (F_a)  \label{eq.linearityMv}
\end{equation}
hold for the matrices $\bM (U_a)$ and $\bM (V_a)$ in Lemma \ref{lemm.sizesuava}. Taking into account \eqref{eq.unsconpoly1}, for each value of $a$, $a = 0,1, \ldots , rd$, we define the function $\Phi \, : \, \mathbb{C}^{p_a \times m} \times \mathbb{C}^{q_a \times n} \rightarrow \mathbb{R}$ as
\begin{equation} \label{eq.phiobjfunc}
\Phi (U_a,V_a) := \left\| 
\left[\begin{array}{c}
P_0 \\
\vdots \\
P_d
\end{array}\right]-
{\bf M}(U_a) \cdot V_a \right\|_F^2 ,
\end{equation}
and consider the optimization problem
\begin{equation}
\label{eq:abstractmin_POLY}
\min_{(U_a ,  V_a )\in \mathbb{C}^{p_a \times m} \times \mathbb{C}^{q_a \times n}} \Phi (U_a,V_a) .
\end{equation}
Observe that we define $\Phi$ as the square of the cost function in \eqref{eq.unsconpoly1}, because we wish to deal with differentiable functions. This does not change the minimizers nor the sequence of matrices $U_a$ and $V_a$ computed by the method in Section \ref{sec:sol_min_pol}. 


Although we deal with complex matrices, we consider $\Phi$ as a function of real variables. More precisely, as a function of the real and imaginary parts of the entries of $U_a$ and $V_a$. In this setting, $\Phi$ is a polynomial of degree four over $\mathbb{R}^{2(p_a m + q_a n)}$ ($\cong \mathbb{C}^{p_a \times m} \times \mathbb{C}^{q_a \times n}$) and, so, it is infinitely differentiable on the whole space. To be consistent with this view, $\mathbb{C}^{p_a \times m} \times \mathbb{C}^{q_a \times n}$ is considered a real vector space endowed with the real inner product
\begin{equation} \label{eq.innerprod}
\langle (A_1,B_1) , (A_2,B_2) \rangle = \Re (\text{trace} (A_1^* A_2)) + 
\Re (\text{trace} (B_1^* B_2))  ,
\end{equation}
where $\Re (z)$ stands for the real part of $z \in \mathbb{C}$.

With the notation and conventions above, the method introduced in Section \ref{sec:sol_min_pol} can be concisely described with the following pseudocode.
\begin{algorithm}[H]
\caption{Alternating least squares for approximating over $\mathcal{C}^{m\times n}_{d,r,a}$} 
\label{alg.conv}
\begin{algorithmic}[1]
\STATE \textbf{Input:} polynomial $P$, tolerance {\it tol}, maximal number of iterations {\it maxiter}
\STATE \textbf{Initialization:} $U_a^{(0)}$
\STATE \textbf{Repeat} recursively generate $(U_a^{(k)}, V_a^{(k)} )$ by 
\medskip
\STATE \qquad  $\displaystyle V_a^{(k)} = \argmin_{V_a \in \mathbb{C}^{q_a \times n}} \Phi(U_a^{(k-1)}, V_a)$
\medskip
\STATE \qquad $\displaystyle U_a^{(k)} = \argmin_{U_a \in \mathbb{C}^{p_a \times m}} \Phi(U_a, V_a^{(k)})$ 
\STATE \textbf{Until} \quad $\left| \sqrt{\Phi(U_a^{(k-1)}, V_a^{(k-1)})} - \sqrt{\Phi(U_a^{(k)}, V_a^{(k)})} \right| \leq tol* \ \left\| P \right\|_F$ or $maxiter$ is reached
\end{algorithmic}
\end{algorithm}

As commented in Section \ref{sec:sol_min_pol}, the problems in lines 4 and 5 of Algorithm \ref{alg.conv} are standard least squares problems and, therefore, have global minimizers for all $k$. Such minimizers are unique if and only if the matrix $\bM (U_a^{(k-1)})$ has full column rank, for the problem in line 4, and if and only if the matrix $\bM (V_a^{(k)})$ has full column rank, for the problem in line~5. Otherwise, there are infinitely many minimizers. Though the generic situation is that the minimizers are unique, as commented in Remark \ref{rem.uniquesizemat}, we do not know if it holds for all $k$, or, more precisely, for which initializations this holds for all $k$. In this setting, we remark that the results in this section are valid even in the case the minimizers are not unique and for any possible choice of them.

Algorithm \ref{alg.conv} is a particular instance of the well-known block coordinate descent methods, whose convergence has been analyzed in depth in the literature. See, for instance, \cite{li2019alternating,xuyin2013-siamimaj}. However, one of the key assumptions in these analyses does not hold for the function $\Phi$ in \eqref{eq.phiobjfunc}, a question that we discuss in Section \ref{subsec.otherconv}. Therefore, we prove directly two results on the convergence of Algorithm \ref{alg.conv}. First, we prove in Theorem \ref{thm.convphi} that the sequence $\left\{ \Phi(U_a^{(k)},V_a^{(k)}) \right\}$ is convergent without the need of any assumption. Second, under the assumption that the sequence  $\left\{ (U_a^{(k)}, V_a^{(k)}) \right\}$ is bounded, we prove in Theorem \ref{thm.convgrad} that the gradient of $\Phi$ converges to zero. Observe that these results do not imply that $\left\{ (U_a^{(k)} , V_a^{(k)} )\right\}$ is convergent. The assumption about the boundedness of $\left\{ ( U_a^{(k)} , V_a^{(k)} )\right\}$ has been imposed in general convergence analyses \cite{li2019alternating,xuyin2013-siamimaj} and has been always satisfied in our numerical tests.
\begin{theorem}
\label{thm.convphi}
Let $\left\{ (U_a^{(k)} , V_a^{(k)} )\right\}$ be the sequence of iterates generated by Algorithm \ref{alg.conv} without stopping criterion and $\left\{ \Phi(U_a^{(k)},V_a^{(k)}) \right\}$ be the sequence of real numbers obtained by evaluating the function $\Phi$ in \eqref{eq.phiobjfunc} at these iterates. Then:
\begin{itemize}
    \item[\rm (a)] $ \Phi(U_a^{(k+1)},V_a^{(k+1)}) \leq  \Phi(U_a^{(k)},V_a^{(k+1)}) \leq  \Phi(U_a^{(k)},V_a^{(k)})$, for $k = 1, 2 , \ldots$.

    \item[\rm (b)] $\left\{ \Phi(U_a^{(k)},V_a^{(k)}) \right\}$  is convergent.
\end{itemize}
\end{theorem}
\begin{proof} Part (a) follows directly from the definition of Algorithm \ref{alg.conv}. Part (a) proves that $\left\{ \Phi(U_a^{(k)},V_a^{(k)}) \right\}$ is a monotone non-increasing sequence. Moreover, it is bounded from below by zero, because
$ \Phi(U^{(k)},V^{(k)}) \geq 0,$ for all $k$, and, therefore, is convergent. 
\end{proof}

\subsection{Convergence to zero of the gradient}
Since the function $\Phi \, : \, \mathbb{C}^{p_a \times m} \times \mathbb{C}^{q_a \times n} \rightarrow \mathbb{R}$ is differentiable everywhere, its differential at $(U_a , V_a)$ is a linear map, denoted by $\mathrm{D}\Phi (U_a , V_a)$, with domain $\mathbb{C}^{p_a \times m} \times \mathbb{C}^{q_a \times n}$ and codomain $\mathbb{R}$. Using the inner product in \eqref{eq.innerprod}, the gradient of $\Phi$ at $(U_a , V_a)$, denoted by $\grad \, \Phi (U_a,V_a)$, is the element of  $\mathbb{C}^{p_a \times m} \times \mathbb{C}^{q_a \times n}$ that satisfies
\begin{equation} \label{eq.defgrad}
\langle \, \grad \, \Phi (U_a,V_a) \, , \, (E_a , F_a) \, \rangle = \mathrm{D}\Phi (U_a , V_a) [(E_a , F_a)] ,
\end{equation}
for all $ (E_a , F_a) \in \mathbb{C}^{p_a \times m} \times \mathbb{C}^{q_a \times n}$.
See \cite[p. 23]{boumalbook-2023} for more details on these concepts. From \eqref{eq.defgrad}, it is possible to obtain the explicit expression of $\grad \, \Phi (U_a,V_a)$.

\begin{lemma} \label{lemm.grad} Let $\Phi \,: \, \mathbb{C}^{p_a \times m} \times \mathbb{C}^{q_a \times n} \rightarrow \mathbb{R}$ be the function in \eqref{eq.phiobjfunc}. Then $\grad \, \Phi (U_a,V_a) = (\, \grad_U \, \Phi (U_a,V_a) \, , \, \grad_V \, \Phi (U_a,V_a) \,)$, with
\begin{align*}
\grad_U \, \Phi (U_a,V_a) & = 2 \,  \bM (V_a)^* \left(  \bM (V_a) \cdot U_a - \begin{bmatrix}
P_0^T \\ \vdots \\ P_d^T
\end{bmatrix} \right) , \\ 
\grad_V \, \Phi (U_a,V_a) & = 2 \, \bM (U_a)^* \left(  \bM (U_a) \cdot V_a - \begin{bmatrix}
P_0 \\ \vdots \\ P_d
\end{bmatrix}\right).
\end{align*}
\end{lemma}
\begin{proof}
We only sketch the proof. It follows from a direct computation of the right-hand side of \eqref{eq.defgrad} as
$$
\mathrm{D}\Phi (U_a , V_a) [(E_a , F_a)] = \lim_{t \rightarrow 0} \frac{\Phi (U_a + t E_a, V_a + t F_a) - \Phi (U_a , V_a )}{t},
$$
taking into account \eqref{eq.phiobjfunc}, that $\| A \|_F^2 = \text{trace} (A^* A)$ for any matrix $A$, and the definition \eqref{eq.innerprod}. In the process, it is needed to apply the equality
\begin{align*}
& \text{trace} \left( \left(  \bM (U_a) \cdot V_a - \begin{bmatrix}
P_0 \\ \vdots \\ P_d
\end{bmatrix}\right)^* \bM(E_a) \cdot V_a \right) = \text{trace} \left( \left(  \bM (V_a) \cdot U_a - \begin{bmatrix}
P_0^T \\ \vdots \\ P_d^T
\end{bmatrix}\right)^* \bM(V_a) \cdot E_a \right),
\end{align*}
which follows from \eqref{eq.reluvvu} and from $\text{trace} (A^* B) = \sum_{i=1}^p  \sum_{j=1}^q \overline{a_{ij}} \, b_{ij}$ for any $A,B \in \mathbb{C}^{p\times q}$.
\end{proof}
Lemma \ref{lemm.grad} proves that $\grad_U \, \Phi (U_a,V_a)$ (resp. $\grad_V \, \Phi (U_a,V_a)$) is two times the difference between the left and the right-hand sides of the normal equations corresponding to the least squares problem \eqref{alt_pol_min2} (resp. \eqref{alt_pol_min1}). Since the normal equations are satisfied by any solution of these least squares problems, we immediately get the following lemma.
\begin{lemma} \label{lemm.grad0pre} Let $\left\{ (U_a^{(k)} , V_a^{(k)} )\right\}$ be the sequence of iterates generated by Algorithm \ref{alg.conv} without stopping criterion. Then
$$
\grad_U \Phi (U_a^{(k)}, V_a^{(k)}) = 0 \qquad \text{and} \qquad \grad_V \Phi (U_a^{(k-1)}, V_a^{(k)}) = 0, \qquad \mbox{for all $k \geq 1$.}
$$
\end{lemma}
Observe that Lemma \ref{lemm.grad0pre} allows us to interpret Algorithm \ref{alg.conv} as a process that alternates between setting to zero $\grad_V \Phi$ and $\grad_U \Phi$.

The rest of this section is devoted to proving, under the assumption that $\left\{ (U_a^{(k)} , V_a^{(k)}) \right\}$ is bounded, that $\lim_{k \rightarrow \infty} \grad_V \Phi (U_a^{(k)}, V_a^{(k)}) = 0$, which, combined with Lemma \ref{lemm.grad0pre}, implies $\lim_{k \rightarrow \infty} \grad  \, \Phi (U_a^{(k)}, V_a^{(k)}) =0$. For that we need the following lemma.
\begin{lemma} \label{lemm.auxxgrado}  Let $\| \cdot \|$ be the norm derived from the inner product \eqref{eq.innerprod} and $\beta$ be a positive number. If the sequence of iterates generated by Algorithm \ref{alg.conv} without stopping criterion satisfies $\| (U_a^{(k)} , V_a^{(k)} )\| \leq \beta$ for all $k \geq 1$, then
\[
\Phi (U_a^{(k)}, V_a^{(k)}) - \Phi (U_a^{(k)}, V_a^{(k+1)}) \geq \frac{1}{4 (d_R +2) \beta^2} \| \grad_V \Phi (U_a^{(k)}, V_a^{(k)}) \|_F^2,
\]
for all $k \geq 1$, where $d_R = \lfloor a/r \rfloor$.  
\end{lemma}
\begin{proof}
Note that $\| (U_a^{(k)} , V_a^{(k)} )\| \leq \beta$ implies that  $\|U_a^{(k)}\|_F^2\leq \beta^2$. Then, Lemma \ref{lemm.sizesuava} and \eqref{eq.matua} imply
$\|\bM(U_a^{(k)})\|_F^2\leq (d_R + 2)\beta^2$, for $k \geq 1$. Let $\bM(U_a^{(k)}) = X \Sigma Y^*$ be an economy-size singular value decomposition of $\bM(U_a^{(k)})$, with $X \in \mathbb{C}^{(d+1)m \times s}$, $\Sigma = \text{diag} (\sigma_1, \ldots , \sigma_s) \in \mathbb{R}^{s \times s}$, $Y \in \mathbb{C}^{q_a\times s}$, and $\sigma_1 \geq \ldots \geq \sigma_s >0$ the nonzero singular values of $\bM(U_a^{(k)})$. Observe that $s = \text{rank} \, (\bM(U_a^{(k)}))$, which can be any value $s \leq (r+a)$. However, if $s=0$, then $\bM(U_a^{(k)})=0$ and $\grad_V \Phi (U_a^{(k)}, V_a^{(k)})=0$, which, combined with Theorem \ref{thm.convphi}-(a), implies trivially the 
result. Therefore, we assume in the rest of the proof that $s \geq 1$. Let $\begin{bmatrix}
X & X_\perp \end{bmatrix} \in \mathbb{C}^{(d+1)m \times (d+1)m}$ be a unitary matrix obtained by completing the columns of $X$. If
$$
\begin{bmatrix}
Z^* \\ Z_\perp^* \end{bmatrix} :=
\begin{bmatrix}
X^* \\ X_\perp^* \end{bmatrix} 
\left[\begin{array}{c}
P_0 \\
\vdots \\
P_d
\end{array}\right],
$$
where $Z^* \in \mathbb{C}^{s \times n}$, then
\begin{align*} \nonumber
\Phi (U_a^{(k)}, V_a^{(k)}) & = 
\left\| \begin{bmatrix}
X^* \\ X_\perp^* \end{bmatrix} \left(
\left[\begin{array}{c}
P_0 \\
\vdots \\
P_d
\end{array}\right]-
{\bf M}(U_a^{(k)}) \cdot V_a^{(k)} \right) \right\|_F^2  = \| Z^* - \Sigma Y^* V_a^{(k)} \|_F^2 + \|Z^*_\perp \|_F^2.
\end{align*}
Analogously,
$$\Phi (U_a^{(k)}, V_a^{(k+1)}) =  \| Z^* - \Sigma Y^* V_a^{(k+1)} \|_F^2 + \|Z^*_\perp \|_F^2 \, . $$
Moreover, in the equality above, $\| Z^* - \Sigma Y^* V_a^{(k+1)}\|_F^2 = 0$, because $V_a^{(k+1)}$ is the output\footnote{Note that $V_a^{(k+1)} = Y \Sigma^{-1} Z^* + W$, with $W$ a matrix whose colums are in the null space of $Y^*$, which is equal to the null space of $\bM(U_a^{(k)})$.} of line 4 in Algorithm \ref{alg.conv} (with $k$ replaced by $k+1$). So,
\[
\Phi (U_a^{(k)}, V_a^{(k)}) - \Phi (U_a^{(k)}, V_a^{(k+1)})  = \| Z^* - \Sigma Y^* V_a^{(k)} \|_F^2 \, .
\]
On the other hand, from Lemma \ref{lemm.grad} and the SVD $\bM(U_a^{(k)}) = X \Sigma Y^*$, we get
\begin{align*}
\| \grad_V \Phi (U_a^{(k)}, V_a^{(k)}) \|_F^2 & = 4 \, \| \Sigma (\Sigma Y^*  V_a^{(k)} - Z^*) \|_F^2 \\
& \leq 4 \, \| \Sigma \|_F^2 \, \|\Sigma Y^*  V_a^{(k)} - Z^*\|_F^2 \\
& = 4 \, \|\bM(U_a^{(k)})\|_F^2 \, (\Phi (U_a^{(k)}, V_a^{(k)}) - \Phi (U_a^{(k)}, V_a^{(k+1)}) ),
\end{align*}
which combined with $\|\bM(U_a^{(k)})\|_F^2\leq (d_R + 2)\beta^2$ completes the proof.
\end{proof} 

A direct corollary of Lemmas \ref{lemm.grad0pre} and \ref{lemm.auxxgrado} and Theorem \ref{thm.convphi} is that the gradient of the cost function converges to zero. This is stated in Theorem \ref{thm.convgrad}.
\begin{theorem} \label{thm.convgrad}  Let $\| \cdot \|$ be the norm derived from the inner product \eqref{eq.innerprod} and $\beta$ be a positive number. If the sequence of iterates generated by Algorithm \ref{alg.conv} without stopping criterion satisfies $\| (U_a^{(k)} , V_a^{(k)} )\| \leq \beta$ for all $k \geq 1$, then
$$
\lim_{k \rightarrow \infty} \grad  \, \Phi (U_a^{(k)}, V_a^{(k)}) =0.
$$
\end{theorem}


\begin{remark} {\rm (On coordinate-wise minima)} If $\{ ( U_a^{(k)} ,\ V_a^{(k)} ) \}$ is bounded, then it has a convergent subsequence  $\{ (U_a^{(k_j)} ,\ V_a^{(k_j)} ) \}$, with limit denoted by  $\{ ( \widetilde U_a , \widetilde V_a ) \}$. Since $\grad \, \Phi$ is a continuous function by Lemma \ref{lemm.grad}, we get from Theorem \ref{thm.convgrad} that $\grad \, \Phi (\widetilde U_a , \widetilde V_a) =0$. The expressions of $\grad_U \Phi (\widetilde U_a , \widetilde V_a)$ and $\grad_V \Phi (\widetilde U_a , \widetilde V_a)$ imply that $\{ ( \widetilde U_a , \widetilde V_a ) \}$ satisfies simultaneously the  normal equations of the least squares problems \eqref{alt_pol_min1} and \eqref{alt_pol_min2}. Thus, $\{ ( \widetilde U_a , \widetilde V_a ) \}$ is a coordinate-wise minimum of $\Phi$ \cite[Definition 14.2]{beckbook-2017}, i.e., $\Phi (\widetilde U_a , \widetilde V_a) \leq \Phi (\widetilde U_a , V_a)$ for any $V_a$ and $\Phi (\widetilde U_a , \widetilde V_a) \leq \Phi (U_a , \widetilde V_a)$ for any $U_a$. Therefore, if $\{ ( U_a^{(k)} ,\ V_a^{(k)} ) \}$ is bounded, any accumulation point of $\{ (U_a^{(k)}  , V_a^{(k)} ) \}$ is a first order stationary point and a coordinate-wise minimum of $\Phi$. 
\end{remark}

\subsection{Relations with general convergence analyses of block coordinate descent methods}
\label{subsec.otherconv} Theorems 2.8 and 2.9 in \cite{xuyin2013-siamimaj} provide results on global convergence and convergence rates of general block coordinate descent methods under certain assumptions and the additional assumption that the sequence of iterates has a finite accumulation point (which is guaranteed by a boundedness assumption as that in Theorem \ref{thm.convgrad}). The three key assumptions on the function $\Phi$ in \eqref{eq.phiobjfunc} 
that would allow us to apply the results in \cite{xuyin2013-siamimaj} are:  (1) the functions $f_k (V_a) := \Phi (U_a^{(k-1)}, V_a)$ and  $g_k (U_a) := \Phi (U_a , V_a^{(k)})$ are strongly convex for all $k$ and their strongly convex parameters have positive lower and upper bounds independent of $k$  \cite[p. 1763]{xuyin2013-siamimaj}; (2) $\grad \, \Phi$ is Lipschitz continuous on any bounded set; and (3) $\Phi$ satisfies the  Kurdyka–\L ojasiewicz property. The function $\Phi$  satisfies (2) because the entries of the gradient are cubic polynomials of the real and imaginary parts of the entries of $U_a$, $V_a$. It also satisfies (3) because $\Phi$ is a quartic polynomial of the real and imaginary parts of the entries of $U_a , V_a$ \cite[p. 1768]{xuyin2013-siamimaj}. However, $\Phi$ does not satisfy (1). Since $f_k (V_a)$ and  $g_k (U_a)$ are quadratic functions, it is not difficult to prove using \cite[Example 5.19]{beckbook-2017} that the condition (1) is equivalent to the existence of two numbers $\ell, L$, $0<\ell < L$, such that the matrices $\bM (U_a^{(k)})$ and $\bM (V_a^{(k)})$ have full column rank and their minimum singular values are included in the interval $[\ell \, , \, L]$ for all $k$. It is easy to construct matrices $\bM (U_a)$ and $\bM (V_a)$ which do not have full column rank. Thus, the required strong convexity cannot be guaranteed in general. However, it might hold in practice, i.e., for most initializations of the algorithm, which would explain the good behavior that we have observed in practice.

\section{Numerical experiments}
\label{sec:numex}

In this section, we compare our algorithm with existing methods.  First, to evaluate if our method can match the best existing approximations, we apply it to matrix pencils and polynomials that have been investigated before. Then we compute the distance to singularity for randomly generated matrix pencils and polynomials using our method and the other state-of-the-art methods.  
The restriction to the distance-to-singularity problem in the comparisons is necessary because mostly the currently available methods do not address the more general problem of approximation by matrix polynomials of arbitrary prescribed rank, such as rank at most $r$. Finally, we illustrate the behavior of our method as the approximating rank $r$ increases. We refer to our method as \textsc{GEAR}, since it is based on GEneric factorizations and Alternating direction method for Rank-$r$ approximation of matrix polynomials. 

All experiments were carried out in MATLAB R2022b Update~8 on a MacBook Pro with an Apple M1 Max processor and 32GB of memory.  The implementation is available on GitHub.\footnote{\url{https://github.com/dmand/GEAR_framework}}

\subsection{Experiments with specific matrix pencils}
 We consider regular matrix pencils $A_0 + \lambda A_1$ taken from different studies on distance to singularity and present the estimated distances for them in Table \ref{tab:pen}. In the last column, we cite the papers where these pencils appear. Our results are typically the same or a bit better than those of the previous studies. Note also that computed distance can often be slightly improved by decreasing the tolerance, but it comes with an increase of the computing time. 

Define the following matrices that will be used as building blocks in the example in row 6 of Table \ref{tab:pen}, 
    \begin{equation*}
    \begin{aligned}
    K_0 &= \begin{bmatrix}
        -67.4894 & -69.2393 & 69.2393\\
        -69.8124 & -1.68624 & 1.68617 \\
        69.8123 & 1.68617 & 68.2707 
    \end{bmatrix}, \quad D_0 = \begin{bmatrix}
        1.52143 & 1.55168 & -1.55168\\
        -3.22064 & -3.28467 & 3.28467  \\
        3.22064 & 3.28467 & -3.28467
    \end{bmatrix}, \\
    M_0 &= \begin{bmatrix}
        18.7532 & -7.94493 & 7.94494\\
        -7.94493 & 31.8182 & -26.8182 \\
        7.94494 & -26.8182 & 26.8182
    \end{bmatrix}, \quad
    F_0 = \begin{bmatrix}
        1 & 0 & 0\\
        0 & 0 & 1
    \end{bmatrix}.
    \end{aligned}
    \end{equation*}
 
\begin{table}[ht!]
\centering
\caption{Computed distances to singularity, $\boldsymbol{d_{\mathrm{sing}}}$, for specific matrix pencils obtained with our method, together with references to previous studies of these problems.}
\label{tab:pen}
\renewcommand{\arraystretch}{0.8}
\setlength{\tabcolsep}{1pt}
\begin{tabular}{c c c c c}
\hline
\textbf{Ex.} & $\boldsymbol{A_0}$ & $\boldsymbol{A_1}$ & $\boldsymbol{d_{\mathrm{sing}}}$ & \textbf{References} \\
\hline
1 &
$\begin{bmatrix}
1 & 1/\varepsilon\\
0 & 1
\end{bmatrix}$ &
$\begin{bmatrix}
0 & 1/\varepsilon\\
0 & 1
\end{bmatrix}$ &
$\varepsilon$ &
\cite{ByHM98, DaBo23, DoNN23}
\\[4mm]

2 &
$\begin{bmatrix}
1 & 0 & 0\\
0 & \varepsilon & 0\\
0 & 0 & 1
\end{bmatrix}$ &
$\begin{bmatrix}
0 & 1 & 0\\
0 & 0 & 1\\
0 & 0 & 0
\end{bmatrix}$ &
$\varepsilon$ &
\cite{ByHM98, DaBo23, DoNN23}
\\[6mm]

3 &
$\begin{bmatrix}
0 & 0.04 & 0.89\\
0.15 & -0.02 & 0\\
0.92 & 0.11 & 0.066
\end{bmatrix}$ &
$\begin{bmatrix}
0 & 0 & 0\\
0 & 0 & 1\\
0 & 1 & 0
\end{bmatrix}$ &
$0.1155462894$ &
\cite{DaBo23, DoNN23, GiHL19, GuLM17}
\\[6mm]

4 &
$\begin{bmatrix}
-1.79 & 0.1 & -0.6\\
0.84 & -0.54 & 0.49\\
-0.89 & 0.3 & 0.74
\end{bmatrix}$ &
$\begin{bmatrix}
0 & 0 & 0\\
0 & 0 & 1\\
0 & 1 & 0
\end{bmatrix}$ &
$0.9435641675$ &
\cite{DaBo23, DoNN23, GiHL19, GuLM17}
\\[6mm]

5 & $\begin{bmatrix}
        1 & -1 & -1 & \dots & -1\\
        & 1 & -1 & \dots & -1\\
        & & \ddots & \ddots & \vdots\\
        & & & 1 & -1\\
         & & & & 1
    \end{bmatrix}$ 
    &
$\begin{bmatrix}
        -1 & 1 & 1 & \dots & 1\\
        & -1 & 1 & \dots & 1\\
        & & \ddots & \ddots & \vdots\\
        & & & -1 & 1\\
         & & & & -1
    \end{bmatrix}$ 
    & $4.0460\times10^{-6}$ 
& \cite{ByHM98, DoNN23}
\\[6mm]
6 & $\begin{bmatrix}
        0 & I_3 & 0\\
        K_0 & D_0 & F_0^T \\
        F_0 & 0 & 0 
    \end{bmatrix}$ 
    &
   $\begin{bmatrix}
        I_3 & 0 & 0 \\
        0 & M_0 & 0\\
        0 & 0 & 0
    \end{bmatrix}$ 
    & $0.0111718453$
    & \cite{ByHM98, DaBo23, DoNN23, GuLM17}
\\[6mm]
\hline
\end{tabular}
\end{table}

\subsection{Experiments with specific matrix polynomials}

In this subsection, we consider several matrix polynomial examples that have previously appeared in the literature. 

\begin{example}\label{ex-7}
Consider the degree--three matrix polynomial 
$ 
P(\lambda)=P_0+\lambda P_1+\lambda^2 P_2+\lambda^3 P_3,
$
where
\[
P_0=\begin{bmatrix}
-0.1414 & -0.149\\
1.1928 & 0.9702
\end{bmatrix},\quad
P_1=\begin{bmatrix}
0.8837 & 0.9969\\
0.219 & 0.0259
\end{bmatrix},
\]
\[
P_2=\begin{bmatrix}
0.6346 & 0.9689\\
0.6252 & -0.0649
\end{bmatrix},\quad
P_3=\begin{bmatrix}
-1.9867 & 1.28\\
0.6097 & -0.1477
\end{bmatrix}.
\]
Using our algorithm, we obtain
$ d_{\mathrm{sing}} = 1.67654181429$ which is consistent with the estimated distance reported in Example~9.3 of \cite{DaBo23}. 
\end{example}

\begin{example}
We consider the quadratic matrix polynomial
$ 
P(\lambda)=P_0+\lambda P_1+\lambda^2 P_2, 
$
with
\[
P_0=\begin{bmatrix}
0.0278 & 0.0563 & 0.1141\\
-0.1758 & 0.327 & -0.173\\
-0.056 & 0.0321 & -0.075
\end{bmatrix},\quad
P_1=\begin{bmatrix}
-0.2122 & 0.363 & -0.1385\\
0.18027 & -0.151 & 0.469\\
-0.106 & 0.212 & -0.1514
\end{bmatrix},
\]
\[
P_2=\begin{bmatrix}
-0.0376 & 0.107 & 0.293\\
0.003 & -0.14914 & -0.2859\\
0.0577 & 0.1455 & 0.231
\end{bmatrix}.
\]
This example was previously studied in Example~9.8 of \cite{DaBo23} and in Section~5.2 of \cite{GiHL17}, where the estimated distance to singularity was reported as $0.026604$, which agrees with our computations. 
\end{example}

\begin{example}\label{best_result}
We consider a degree--three matrix polynomial
$
P(\lambda)=P_0+\lambda P_1+\lambda^2 P_2+\lambda^3 P_3, 
$
where 
    $$\small{P_3 = \begin{bmatrix}
        0 & 0.01837967 & 0 & 0\\
        0 & 0.01843603 & 0 & 0\\
        0 & 0.01829203 & 0 & 0\\
        0 & 0.01842778 & 0 &0
    \end{bmatrix}, \ \ \small{P_2 = \begin{bmatrix}
        0.01833149 & 0.0366177 & 0.01824331 & 0.03660918\\
        0.01837542 & -0.05442525 & 0 & 0.01832234\\
        0.01841784 & 0.00003900436 & 0 & 0.0183651578\\
        0.01840752 & 0.00001508311 & 0.01839699 & 0.0365917
    \end{bmatrix}},}
    $$ 
    $$\small{P_1 = \begin{bmatrix}
        0.09116729 & 0.0000179769 & 0.2550857 & 0.05475106\\
        0.0001156514 & 0.00001659159 & 0.09108906 & -0.05447104\\
        0.05470823 & 0.03662426 & 0.1276959 & 0.03650378\\
        0.05472202 & -0.1091389 & 0.1458359 & -0.09090507
    \end{bmatrix},}$$
    $$
    \small{P_0 = \begin{bmatrix}
        0.09108776 & -0.05442464 & 0.3645006 & 0.01821543\\
        -0.1456436 & 0.03647524 & -0.07277662 & 0.07305016\\
        0.05478714 & -0.05444916 & 0.437322 & 0.05478385\\
        -0.1274211 & 0.09124859 & -0.6556615 &-0.0544685
    \end{bmatrix}.}$$ 
The distances to rank-2 approximations reported in the literature for this polynomial are $7.7049\times10^{-4}$ in \cite[Example~9.1]{DaBo23} and $7.844\times10^{-4}$ in \cite[Example~5.6]{GiHL19}. Once again, our algorithm matches the best reported result, computing an approximation
at a distance of $7.7049\times10^{-4}$.
\end{example}

\subsection{Speed and quality for random matrix pencils}

\begin{example}
\label{ex:small}
We calculate a distance to singularity using the three potentially fastest available algorithms, i.e., GEAR (our algorithm), Riemannian \cite{DoNN23}, and BFGS \cite{DaBo23}, for 40 random matrix pencils of each size (20 real and 20 complex), with the sizes from $5 \times 5$ to $35 \times 35$. The results for matrix pencils with real entries are summarized in Figure \ref{fig:5-35-real} and Table \ref{tab:5-35-real}. The graph on the left in Figure \ref{fig:5-35-real} shows that for pencils of these sizes our algorithm is significantly (100 times) faster than competitors. Moreover, the graph on the right in Figure \ref{fig:5-35-real} and Table \ref{tab:5-35-real} show that the computed distance to singularity by our algorithm is typically better (smaller).

\begin{figure}[ht]
    \hspace{-1.7cm}\includegraphics[width=1.13\textwidth]{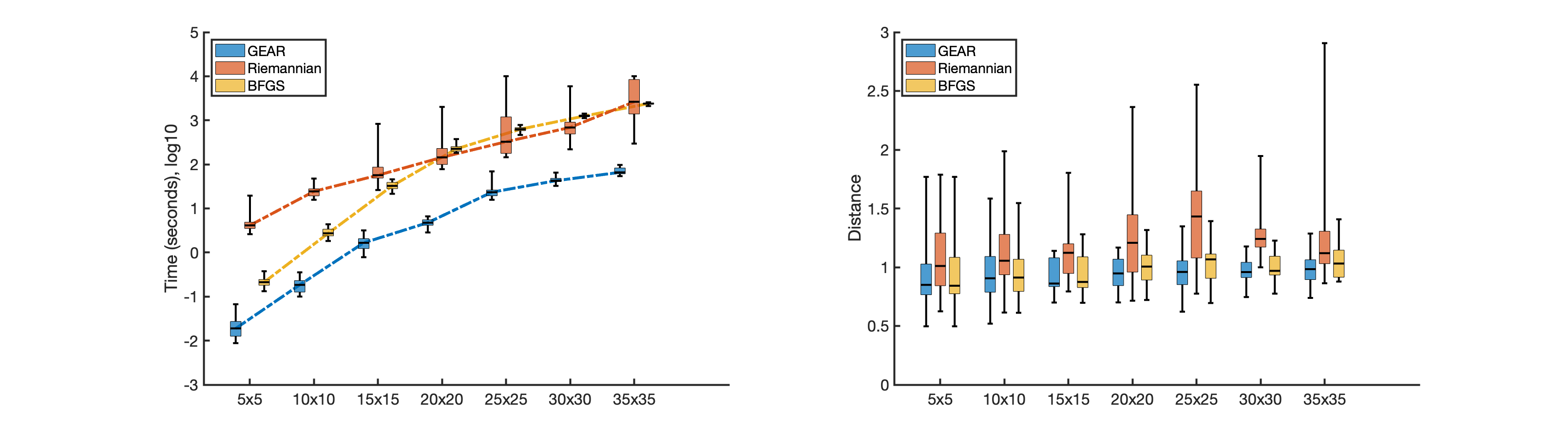}
    \caption{Comparison of the speed and computed minimal distances for random real matrix pencils of given size.}
    \label{fig:5-35-real}
\end{figure}

\begin{table}[ht]
  \centering
  \renewcommand{\arraystretch}{1.2} 
  \scalemath{0.92}{
  \begin{tabular}{l|c|c|c|c|c|c|c}
    \toprule
     & $5 \times 5$ & $10 \times 10$ & $15 \times 15$ & $20 \times 20$ & $25 \times 25$ & $30 \times 30$ & $35 \times 35$ \\
    \midrule
    GEAR         & 4            & 10             & 10             & 17             & 15             & 14             & 18             \\
    Riemannian & 1            & 0              & 1              & 2              & 1              & 1              & 1              \\
    BFGS       & 15           & 10             & 9              & 1              & 4              & 5              & 1              \\
    \bottomrule
  \end{tabular}}
  \smallskip
  \caption{The table shows how many times each algorithm, i.e., GEAR, Riemannian, or BFGS (listed in the first column), returned the best answer out of 20 experiments for random real matrix pencils of a corresponding size (between 5 and 35, all sizes are listed in the first row).}
  \label{tab:5-35-real}
\end{table}

The results for matrix pencils with complex entries are summarized in Figure \ref{fig:5-35-complex} and Table \ref{tab:5-35-complex}. The graph on the left in Figure \ref{fig:5-35-complex} shows that for pencils of sizes up to $30 \times 30$ our algorithm is faster than competitors. For $35 \times 35$ pencils Riemannian is the fastest algorithm, but the variation of its outputs is very large comparing to the other two algorithms. The graph on the right in Figure \ref{fig:5-35-complex} and Table \ref{tab:5-35-complex} show that the computed distance to singularity by our algorithm is typically better (smaller).  

\begin{figure}[ht]
    \hspace{-1.7cm}\includegraphics[width=1.115\textwidth]{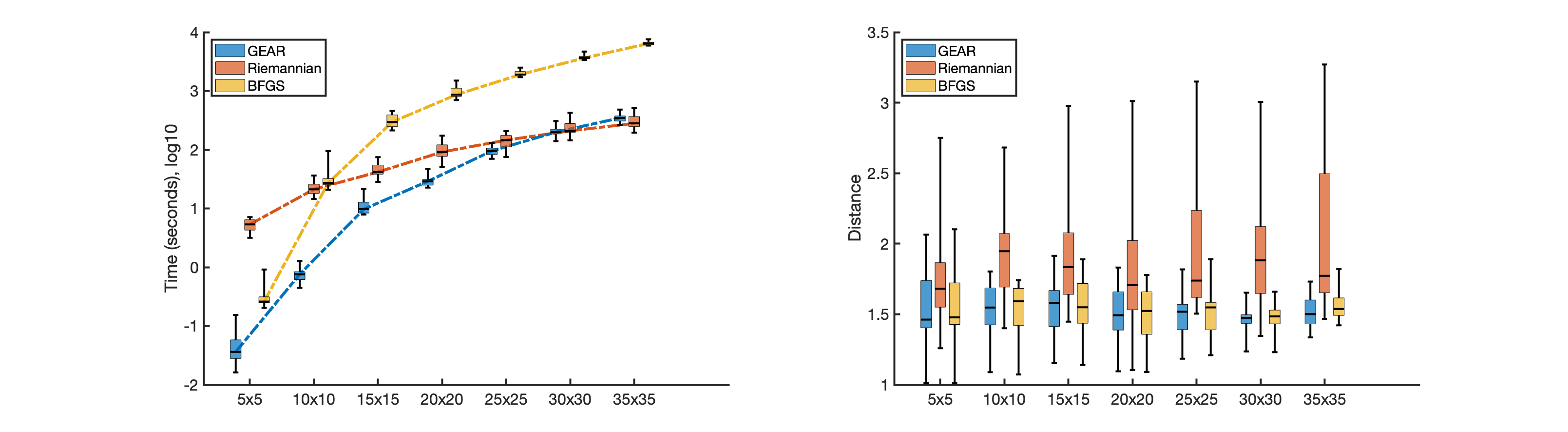}
    \caption{Comparison of the speed and computed minimal distances for random complex matrix pencils of given size.}
    \label{fig:5-35-complex}
\end{figure}

\begin{table}[ht]
  \centering

  \renewcommand{\arraystretch}{1.0} 
\scalemath{0.92}{
  \begin{tabular}{l|c|c|c|c|c|c|c}
    \toprule
     & $5 \times 5$ & $10 \times 10$ & $15 \times 15$ & $20 \times 20$ & $25 \times 25$ & $30 \times 30$ & $35 \times 35$ \\
    \midrule
    GEAR   &    3   &  7   &  7  &   8  &  11  &  13  &  18 \\
    Riemannian &0  &   0  &   1   &  0   &  0 &    0  &   1  \\
    BFGS  & 17  &  13 &   12  &  12  &   9  &   7  &   1 \\
    \bottomrule
  \end{tabular}}
  \smallskip
  \caption{The table shows how many times each algorithm, i.e., GEAR, Riemannian, or BFGS (listed in the first column), returned the best answer out of 20 experiments for random complex matrix pencils of a corresponding size (between 5 and 35, all sizes are listed in the first row).}
  \label{tab:5-35-complex}
\end{table}
\end{example}

\begin{example}
We calculate a distance to singularity using the two fastest algorithms (based on Example \ref{ex:small}), i.e., GEAR (our algorithm) and Riemannian \cite{DoNN23}, for 20 random matrix pencils of each size (10 real and 10 complex), with the sizes $50 \times 50$, $75 \times 75$, and $100 \times 100$. The results for matrix pencils with real and complex entries are summarized in Figures \ref{fig:50-75-100-real} and \ref{fig:50-75-100-complex}, respectively. 
The plots in Figure~\ref{fig:50-75-100-real} indicate that, for real pencils of these sizes, GEAR is consistently faster than Riemannian and we checked that each time it computes a smaller distance to singularity.
Figure \ref{fig:50-75-100-complex} shows that for complex pencils of these sizes our algorithm with the tolerance parameter set to $10^{-7}$ is slower than Riemannian but consistently produces approximations at a smaller distance to singularity, while our algorithm with the tolerance parameter set to $10^{-5}$ is slightly faster than Riemannian; however, as the pencil size increases, it more frequently produces approximations at a larger distance than those obtained by Riemannian. 
\begin{figure}[ht]
\hspace{-1.7cm}\includegraphics[width=1.115\textwidth]{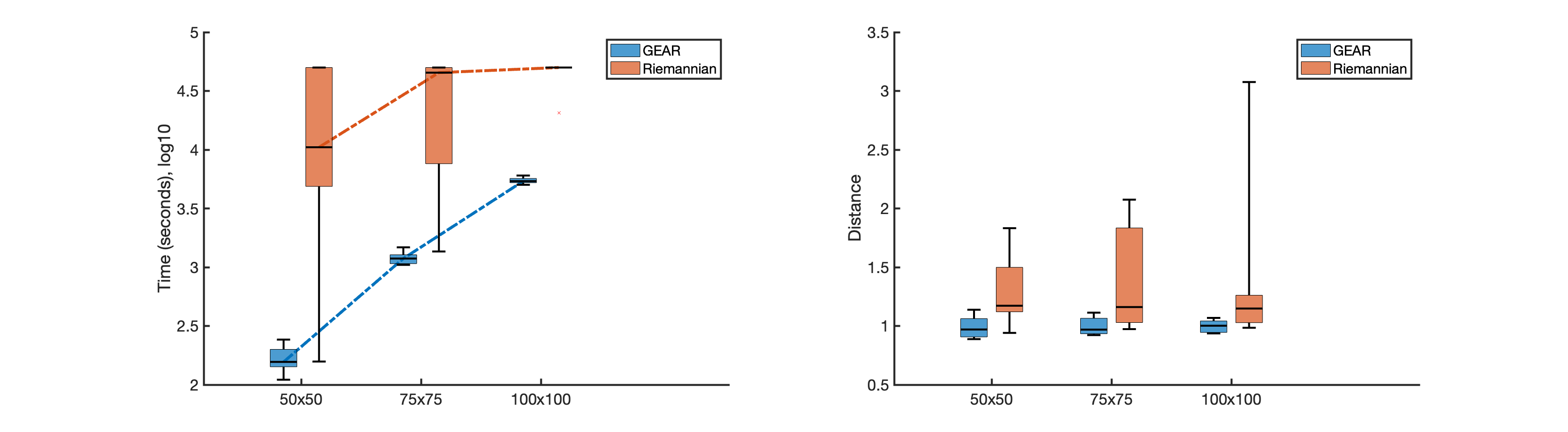}
    \caption{Comparison of the speed and computed minimal distances for random real matrix pencils of given size. For a large fraction of experiments and, in particular for all the $100 \times 100$ experiments, Riemannian reached a maximum allowed time without converging (increasing the maximum allowed time did not improve the situation). In addition, we checked that GEAR found a smaller minimum than Riemannian in each of the 30 performed experiments.
    }
    \label{fig:50-75-100-real}
\end{figure}
\begin{figure}[ht]
    
\hspace{-1.7cm}\includegraphics[width=1.115\textwidth]{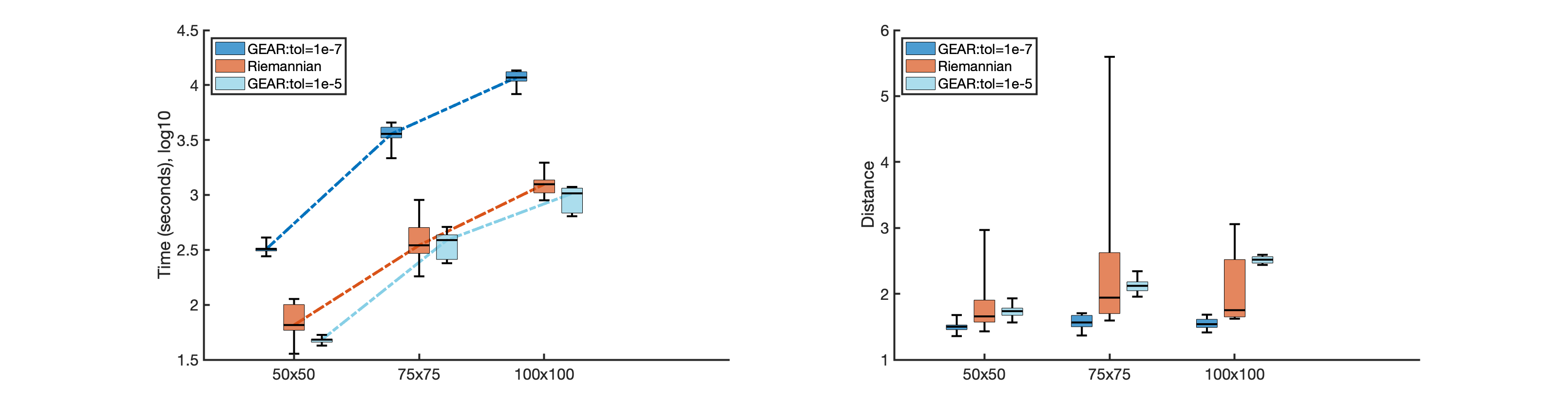}
    \caption{Comparison of the speed and computed minimal distances for random complex matrix pencils of given size. In addition, we checked that GEAR with the tolerance parameter set to $10^{-7}$ produced approximations at the smallest distance to singularity in each of the 30 performed experiments.}
    \label{fig:50-75-100-complex}
\end{figure}
\end{example}

\subsection{Speed and quality for random matrix polynomials}

\begin{example}
\label{ex:pol_small}
We calculate a distance to singularity using four available algorithms, i.e., GEAR (our algorithm), MVF \cite{GnGu23}, Riemann-Oracle (RO) \cite{GNNP24}, and BFGS \cite{DaBo23}, for 40 random matrix polynomials of degree $3$, and of sizes from $3 \times 3$ to $15 \times 15$ (we take 16 polynomials of each size, 8 real and 8 complex). The results for matrix polynomials with real entries are summarized in Figure \ref{fig:4-12-deg3}. The graph on the left in Figure \ref{fig:4-12-deg3} shows that for polynomials of these sizes our algorithm is the fastest among the tested algorithms.
To compare the quality of the minima obtained by the algorithms, we now assess the ``singularity'' of the obtained output, as seen in the graph on the right in Figure \ref{fig:4-12-deg3}, as well as the distance to it, which is shown in the middle graph in Figure \ref{fig:4-12-deg3} (comparing just the latter was sufficient in the pencil case). Our algorithm and BFGS produce the answers that are ``the most singular''.

\begin{figure}[ht]
    \hspace{-1.7cm}\includegraphics[width=1.115\textwidth]{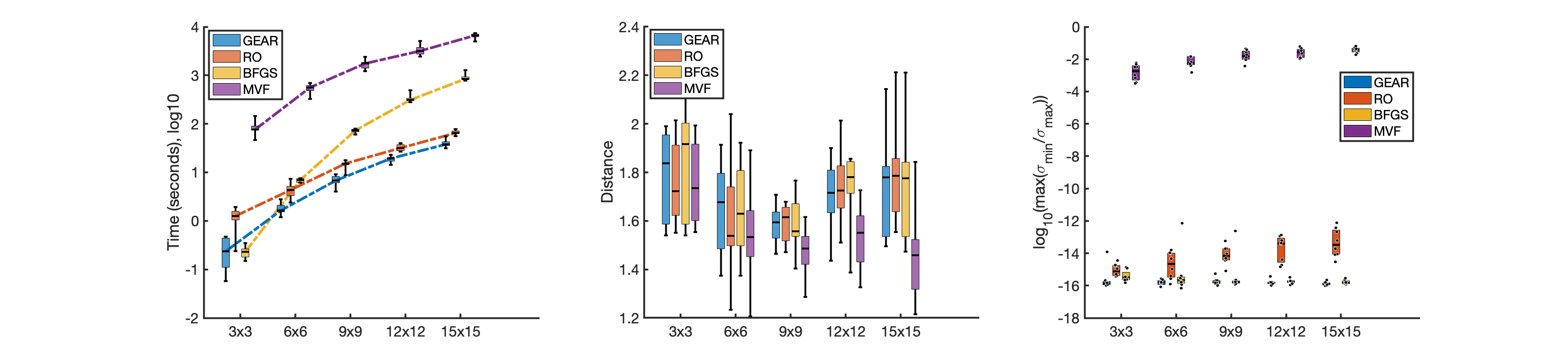}
    \caption{Comparison of the speed, computed minimal distances to singularity, and singularity of the outputs for random real matrix polynomials of degree 3 and given size. 
    To compare the ``singularity'' of the output, we plot, in the right-most graph, the maximum value of the ratio of the smallest singular value $\sigma_{\min}$ to the largest singular value $\sigma_{\max}$ of the constructed matrix polynomial evaluated at the points $x_j+\sqrt{-1}\,y_j$, where $[X,Y]=\texttt{meshgrid}(-100:20:100)$. In all the experiments performed, GEAR produces outputs with this ratio at the level of machine precision.
    }
    \label{fig:4-12-deg3}
\end{figure}

The results for matrix polynomials with complex entries are summarized in Figure~\ref{fig:3-12-complex}. The left plot in Figure~\ref{fig:3-12-complex} shows that, for polynomials of these sizes, our algorithm is the fastest. To compare the quality of the minima obtained by the different algorithms, we again assess both the ``singularity'' of the computed outputs, shown in the right plot of Figure~\ref{fig:3-12-complex}, and their distance to the original polynomial, shown in the middle plot. Our algorithm and BFGS perform best in preserving the singularity, followed by RO.
The latter three algorithms produce outputs at approximately the same distance from the original polynomial.

\begin{figure}[ht]
    \hspace{-1.7cm}\includegraphics[width=1.115\textwidth]{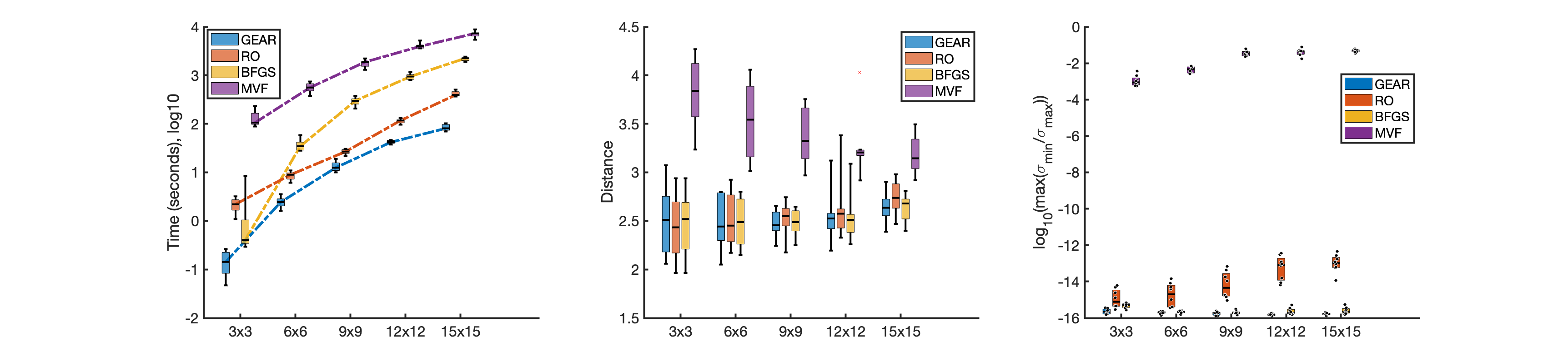}
    \caption{Comparison of the speed, computed minimal distances to singularity, and singularity of the outputs for random complex matrix polynomials of degree 3 and given size. To compare the ``singularity'' of the output, we plot, in the right-most graph, the maximum value of the ratio of the smallest singular value $\sigma_{\min}$ to the largest singular value $\sigma_{\max}$ of the constructed matrix polynomial evaluated at the points $x_j+\sqrt{-1}\,y_j$, where $[X,Y]=\texttt{meshgrid}(-100:20:100)$. In all the experiments performed, GEAR produces outputs with this ratio at the level of machine precision.}
    \label{fig:3-12-complex}
\end{figure}

\end{example}

\begin{example}
\label{ex:ranksvar}
We illustrate how the time needed for computing an approximation with GEAR decreases with the decrease of the rank. 
As an example, we use 60 random real matrix polynomials of degree 2 and sizes $70 \times 70$, $120 \times 120$, and $170 \times 170$ (4 matrix polynomials per each size and rank). In particular, in all the cases considered, computing a rank-5 approximation is around 100 times faster than computing a rank-45 approximation, see more in Figure \ref{fig:GEAR_ranks}. 
    \begin{figure}[ht]  \hspace{-1.2cm}\includegraphics[width=1.0\textwidth]{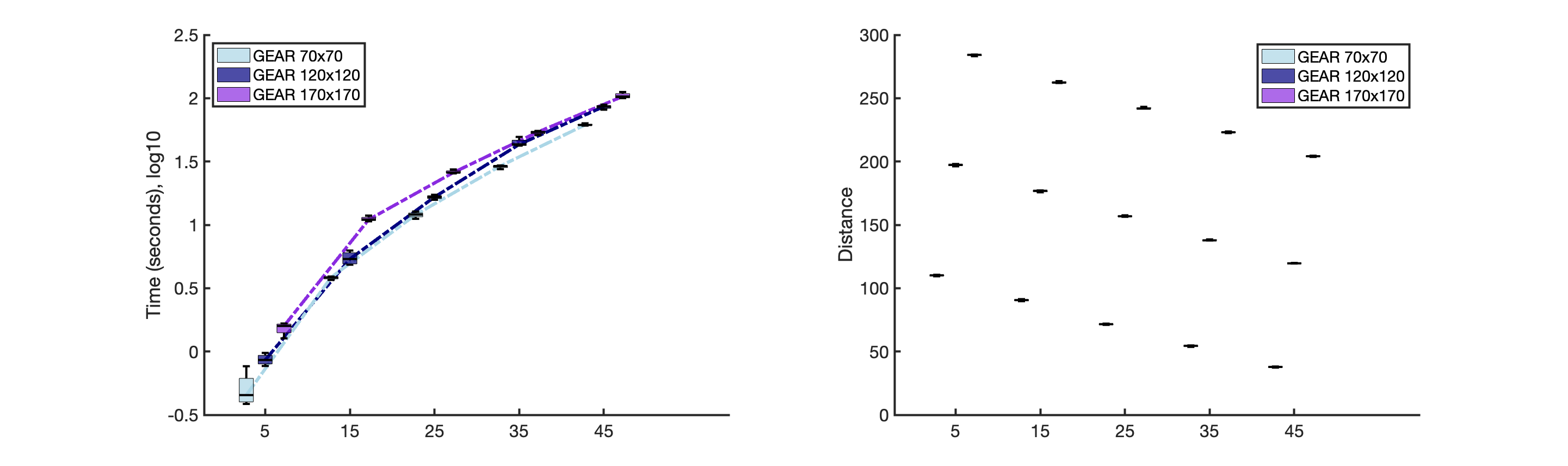}
    \caption{Speed and computed minimal distances to polynomials of fixed rank (varying from 5 to 45) for random real matrix polynomials of degree 2 and sizes $70 \times 70, 120 \times 120$, and $170 \times 170$.}
    \label{fig:GEAR_ranks}
\end{figure}
\end{example}

\section{Conclusions and future work}
\label{sec:conclusions}

We have introduced a novel framework for low-rank approximation of matrix polynomials.  By using explicit parametrizations of generic orbits of matrix polynomials with prescribed normal rank, the original constrained problem is transformed into an unconstrained nonlinear least-squares problem. One of the main advantages of this formulation is that the rank constraint is enforced by construction: every polynomial generated by the method has normal rank at most the prescribed target rank.  The framework applies to arbitrary degrees and to arbitrary ranks $r<\min(m,n)$, and is not restricted to the distance-to-singularity case. The alternating least-squares method in this paper provides a simple and practical realization of this idea, and the numerical experiments show that it is very effective. 

The proposed framework is modular.  Different generic rank-$r$ families, different parametrizations, and different solvers can be combined following the same approach. Therefore, our framework opens up multiple directions for future research. The most pressing one is the development of faster and more scalable optimization methods.  Speed is still a bottleneck for large problems; improving it would make it possible to handle larger matrix polynomials and to use multiple starting guesses for problems of moderate size. Since the resulting nonlinear least-squares problems are nonconvex, multiple initializations may be important for obtaining better approximations.
Other important directions include the design of reliable initialization strategies, the development of structure-preserving methods (feasibility and efficiency of this direction has been showed in a recent preprint \cite{DmDH26}), and the study of alternative parametrizations. Moreover, improving our understanding of the geometry of the involved objects may also lead to sharper convergence results, improved stopping criteria, and better indicators of the quality of the computed approximations.




\begin{funding}
S. Das and A. Dmytryshyn were supported by the Swedish Research Council (VR) grant 2021-05393.
F. M. Dopico was funded by MICIU/AEI/10.13039/501100011033/ and by FEDER/UE through grant PID2023-147366NB-I00.
\end{funding}


{\footnotesize
\bibliographystyle{abbrv}
\bibliography{library}

@article{DmDH26,
  title={Distance to nearest skew-symmetric matrix polynomials of bounded rank},
  author={Dmytryshyn, Andrii and Dopico, Froil{\'a}n M and Hellberg, Rakel},
  journal={arXiv preprint arXiv:2601.16676},
  year={2026}
}

@article{GIHL19,
  title={Computing lower rank approximations of matrix polynomials},
  author={Giesbrecht, Mark and Haraldson, Joseph and Labahn, George},
  journal={Journal of Symbolic Computation},
  volume={98},
  pages={225--245},
  year={2020},
  publisher={Elsevier}
}

@inproceedings{SrJa03,
  title={Weighted low-rank approximations},
  author={Srebro, Nathan and Jaakkola, Tommi},
  booktitle={Proceedings of the 20th International Conference on Machine Learning (ICML-03)},
  pages={720--727},
  year={2003}
}

@InProceedings{LiLR16,
  title = 	 {Recovery guarantee of weighted low-rank approximation via alternating minimization},
  author = 	 {Li, Yuanzhi and Liang, Yingyu and Risteski, Andrej},
  booktitle = 	 {Proceedings of The 33rd International Conference on Machine Learning},
  pages = 	 {2358--2367},
  year = 	 {2016},
  editor = 	 {Balcan, Maria Florina and Weinberger, Kilian Q.},
  volume = 	 {48},
  series = 	 {Proceedings of Machine Learning Research},
  address = 	 {New York, New York, USA},
  month = 	 {20--22 Jun},
  publisher =    {PMLR},
  url = 	 {https://proceedings.mlr.press/v48/lii16.html}
}

@book{beckbook-2017,
 author = {Beck, Amir},
 title = {First-order methods in optimization},
 fseries = {MOS-SIAM Series on Optimization},
 series = {MOS-SIAM Ser. Optim.},
 volume = {25},
 isbn = {978-1-61197-498-0; 978-1-61197-499-7},
 year = {2017},
 publisher = {Philadelphia, PA: Society for Industrial {and} Applied Mathematics (SIAM)},
 language = {English},
 doi = {10.1137/1.9781611974997},
 zbMATH = {6821287},
 Zbl = {1384.65033}
}

@book{boumalbook-2023,
 author = {Boumal, Nicolas},
 title = {An introduction to optimization on smooth manifolds},
 isbn = {978-1-00-916617-1; 978-1-00-916615-7; 978-1-00-916616-4},
 year = {2023},
 publisher = {Cambridge: Cambridge University Press},
 language = {English},
 doi = {10.1017/9781009166164},
 zbMATH = {7633911},
 Zbl = {1532.90001}
}

@inproceedings{li2019alternating,
  title={Alternating minimizations converge to second-order optimal solutions},
  author={Li, Qiuwei and Zhu, Zhihui and Tang, Gongguo},
  booktitle={International Conference on Machine Learning},
  pages={3935--3943},
  year={2019},
  organization={PMLR}
}

@article{xuyin2013-siamimaj,
 author = {Xu, Yangyang and Yin, Wotao},
 title = {A block coordinate descent method for regularized multiconvex optimization with applications to nonnegative tensor factorization and completion},
 fjournal = {SIAM Journal on Imaging Sciences},
 journal = {SIAM J. Imaging Sci.},
 issn = {1936-4954},
 volume = {6},
 number = {3},
 pages = {1758--1789},
 year = {2013},
 language = {English},
 doi = {10.1137/120887795},
 zbMATH = {6249711},
 Zbl = {1280.49042}
}

@article{GnGu25,
  title={On the numerical approximation of the distance to singularity for matrix-valued functions},
  author={Gnazzo, Miryam and Guglielmi, Nicola},
  journal={SIAM J. Matrix Anal. Appl.},
  volume={46},
  number={2},
  pages={1484--1517},
  year={2025},
  publisher={SIAM}
}

@article{GNNP24,
  title={Riemann-Oracle: A general-purpose {R}iemannian optimizer to solve nearness problems in matrix theory},
  author={Gnazzo, Miryam and Noferini, Vanni and Nyman, Lauri and Poloni, Federico},
  journal={Found. Comput. Math.},
  pages={1--45},
  year={2025},
  publisher={Springer}
}

@article{DmDV23,
    AUTHOR = {Dmytryshyn, Andrii and Dopico, Froil\'an M. and Van Dooren,
              Paul},
     TITLE = {Minimal rank factorizations of polynomial matrices},
   JOURNAL = {Linear Algebra Appl.},
  FJOURNAL = {Linear Algebra and its Applications},
    VOLUME = {721},
      YEAR = {2025},
     PAGES = {736--774},
      ISSN = {0024-3795,1873-1856},
   MRCLASS = {15A54 (15A18 15A22 15A23)},
  MRNUMBER = {4913758},
       DOI = {10.1016/j.laa.2025.01.009},
       URL = {https://doi.org/10.1016/j.laa.2025.01.009},
}

@incollection{KrVo15,
    AUTHOR = {Kressner, Daniel and Voigt, Matthias},
     TITLE = {Distance problems for linear dynamical systems},
 BOOKTITLE = {Numerical {A}lgebra, {M}atrix {T}heory, {D}ifferential-{A}lgebraic
              {E}quations and {C}ontrol {T}heory},
     PAGES = {559--583},
 PUBLISHER = {Springer, Cham},
      YEAR = {2015},
      ISBN = {978-3-319-15259-2; 978-3-319-15260-8},
   MRCLASS = {93C73 (93D09)},
  MRNUMBER = {3380294},
MRREVIEWER = {Thomas\ Berger},
}

@article{GnGu23,
      title={Approximating the closest structured singular matrix polynomial}, 
      author={Miryam Gnazzo and Nicola Guglielmi},
      journal = {Linear and Multilinear Algebra},
      volume = {73},
      number = {9},
      pages = {2222--2250},
      year = {2025} 
}

@article{DoNN23,
      title={A {R}iemannian optimization method to compute the nearest singular pencil}, 
      author={Froilán M. Dopico and Vanni Noferini and Lauri Nyman},
 journal={SIAM J. Matrix Anal. Appl.},
volume = {45},
number = {4},
pages = {2007-2038},
  year={2024},
  publisher={SIAM}

}

@article{DaBo23,
title = {Nearest rank deficient matrix polynomials},
journal = {{L}inear {A}lgebra {A}ppl.},
volume = {674},
pages = {304-350},
year = {2023},
issn = {0024-3795},
doi = {https://doi.org/10.1016/j.laa.2023.05.019},
url = {https://www.sciencedirect.com/science/article/pii/S0024379523002021},
author = {Das, Biswajit  and Bora, Shreemayee },
}

@inproceedings{GiHL17,
 author = {Giesbrecht, Mark and Haraldson, Joseph and Labahn, George},
 title = {Computing the Nearest Rank-Deficient Matrix Polynomial},
 booktitle = {Proceedings of the 2017 ACM International Symposium on Symbolic and Algebraic Computation},
 series = {ISSAC '17},
 year = {2017},
 isbn = {978-1-4503-5064-8},
 location = {Kaiserslautern, Germany},
 pages = {181--188},
 numpages = {8},
 url = {http://doi.acm.org/10.1145/3087604.3087648},
 doi = {10.1145/3087604.3087648},
 acmid = {3087648},
 publisher = {ACM},
 address = {New York, NY, USA},
}

@article{GuLM17,
  title={On the Nearest Singular Matrix Pencil},
  author={Nicola Guglielmi and Christian Lubich and Volker Mehrmann},
  journal={SIAM J. Matrix Anal. Appl.},
  year={2017},
  volume={38},
  pages={776--806}
}

@article{ByHM98,
title = "Where is the nearest non-regular pencil?",
journal = "{L}inear {A}lgebra {A}ppl.",
volume = "285",
number = "1",
pages = "81-105",
year = "1998",
issn = "0024-3795",
doi = "https://doi.org/10.1016/S0024-3795(98)10122-2",
url = "http://www.sciencedirect.com/science/article/pii/S0024379598101222",
author = "Ralph Byers and Chunyang He and Volker Mehrmann"
}

@Article{DmDo17,
  author      = {Dmytryshyn, A. and Dopico, F. M.},
  title       = {{G}eneric complete eigenstructures for sets of matrix polynomials with bounded rank and degree},
  year        = {2017},
  journal   = {{L}inear {A}lgebra {A}ppl.},
  volume = {535},
  year = {2017},
  pages = {213--230},
  doi = {10.1016/j.laa.2017.09.007},
}

@article{DeDL17,
author = {De Ter\'an, F. and Dopico, F. M. and Landsberg, J. M.},
title = {An explicit description of the irreducible components of the set of matrix pencils with bounded normal rank},
journal = {Linear Algebra Appl.},
volume = {520},
year = {2017},
pages = {80--103},
issn = {0024-3795},
}

@article {DeDo16,
    AUTHOR = {De Ter{\'a}n, Fernando and Dopico, F. M.},
     TITLE = {Generic change of the partial multiplicities of regular matrix
              pencils under low-rank perturbations},
   JOURNAL = {SIAM J. Matrix Anal. Appl.},
  FJOURNAL = {SIAM Journal on Matrix Analysis and Applications},
    VOLUME = {37},
      YEAR = {2016},
    NUMBER = {3},
     PAGES = {823--835},
      ISSN = {0895-4798},
     CODEN = {SJMAEL},
   MRCLASS = {15A22 (15A18 15A21 65F15)},
  MRNUMBER = {3521544},
       DOI = {10.1137/15M1022069},
       URL = {http://dx.doi.org/10.1137/15M1022069},
}

@Article{BHMS13,
  author     = {Betcke, T. and Higham, N.J. and Mehrmann, V. and Schr\"{o}der, C. and Tisseur, F.},
  title      = {{NLEVP}: {A} Collection of Nonlinear Eigenvalue Problems},
  year       = {2013},
  volume     = {39},
  number     = {2},
  pages      = {7:1--7:28},
  issn       = {0098-3500},
  doi        = {10.1145/2427023.2427024},
  acmid      = {2427024},
  address    = {{N}ew {Y}ork, {NY}, {USA}},
  articleno  = {7},
  issue_date = {February 2013},
  journal    = {{ACM} {T}rans. {M}ath. {S}oftware},
  numpages   = {28},
  publisher  = {ACM},
}

@Book{GoVa13,
  author    = {Golub, G.H. and {Van~Loan}, C.F.},
  title     = {{M}atrix {C}omputations},
  year      = {2013},
  edition   = {4rd},
  publisher = {John Hopkins University Press},
  address   = {{B}altimore, {MD}},
}

@Article{TiMe01,
  author    = {Tisseur, F. and Meerbergen, K.},
  title     = {{T}he quadratic eigenvalue problem},
  year      = {2001},
  volume    = {43},
  number    = {2},
  pages     = {235--286},
  journal   = {{SIAM} {R}ev.},
}

@InCollection{MaMT15,
  author    = {Mackey, D. S. and Mackey, N. and Tisseur, F.},
  title     = {{P}olynomial {E}igenvalue {P}roblems: {T}heory, {C}omputation, and {S}tructure},
  booktitle = {{N}umerical {A}lgebra, {M}atrix {T}heory, {D}ifferential-{A}lgebraic {E}quations and {C}ontrol {T}heory},
  year      = {2015},
  publisher = {Springer, Cham},
  pages     = {319--348},
}
}

\end{document}